\UseRawInputEncoding

\documentclass[final]{elsarticle}
\usepackage{amsmath}
\usepackage{float}
\usepackage{amsfonts}
\usepackage{bm}
\usepackage{amssymb}
\usepackage{chngpage}
\usepackage{multirow}
\usepackage{natbib}
\setcitestyle{comma,numbers,super,open={},close={}}
\usepackage{geometry}
\usepackage[linesnumbered, ruled,vlined]{algorithm2e}
\usepackage[export]{adjustbox}
\usepackage{graphicx}
\usepackage{epstopdf}
\usepackage{rotating}
\usepackage{lineno}
\usepackage{hyperref}
\usepackage{makecell}

\SetKwInput{KwInput}{Input}                
\SetKwInput{KwOutput}{Output}              
\SetKwInput{KwInitialise}{Initialisation}  

\usepackage[table,xcdraw]{xcolor}
\usepackage{graphicx}
\usepackage{subcaption}
\usepackage{mwe}
\usepackage{setspace}
\usepackage{blindtext}
\usepackage{bbm}

\hypersetup{colorlinks = true,linkcolor = blue,anchorcolor =red,citecolor = blue,filecolor = red,urlcolor = red,  pdfauthor=author}

\journal{XXX} 

\begin{document}

\begin{frontmatter}

\title{Hierarchical planning-scheduling-control - Optimality surrogates and derivative-free optimization}

\author[a]{Damien van de Berg}
\author[a]{Nilay Shah\corref{coraut1}}
 \ead{n.shah@imperial.ac.uk}
 \cortext[coraut1]{Corresponding Author}
\author[a]{Ehecatl Antonio del Rio-Chanona\corref{coraut2}}
\ead{a.del-rio-chanona@imperial.ac.uk}
\cortext[coraut2]{Corresponding Author}

\address[a]{Sargent Centre for Process Systems Engineering, Roderic Hill Building
South Kensington Campus
London, SW7 2AZ, United Kingdom}


\doublespacing
\small

%
\begin{abstract}
Planning, scheduling, and control typically constitute separate decision-making units within chemical companies. Traditionally, their integration is modelled sequentially, but recent efforts prioritize lower-level feasibility and optimality, leading to large-scale, potentially multi-level, hierarchical formulations. Data-driven techniques, like optimality surrogates or derivative-free optimization, become essential in addressing ensuing tractability challenges. We demonstrate a step-by-step workflow to find a tractable solution to a tri-level formulation of a multi-site, multi-product planning-scheduling-control case study. We discuss solution tractability-accuracy trade-offs and scaling properties for both methods. Despite individual improvements over conventional heuristics, both approaches present drawbacks. Consequently, we synthesize our findings into a methodology combining their strengths. Our approach remains agnostic to the level-specific formulations when the linking variables are identified and retains the heuristic sequential solution as fallback option. We advance the field by leveraging parallelization, hyperparameter tuning, and a combination of off- and on-line computation, to find tractable solutions to more accurate \textit{multi-level} formulations.

\end{abstract}

\end{frontmatter}

\doublespacing

\nolinenumbers
\small

\makeatletter
\@fpsep\textheight
\makeatother

\section{Introduction}

\subsection{Background}
\begin{figure}[htp!] 
    \centering
    \includegraphics[width=\textwidth]{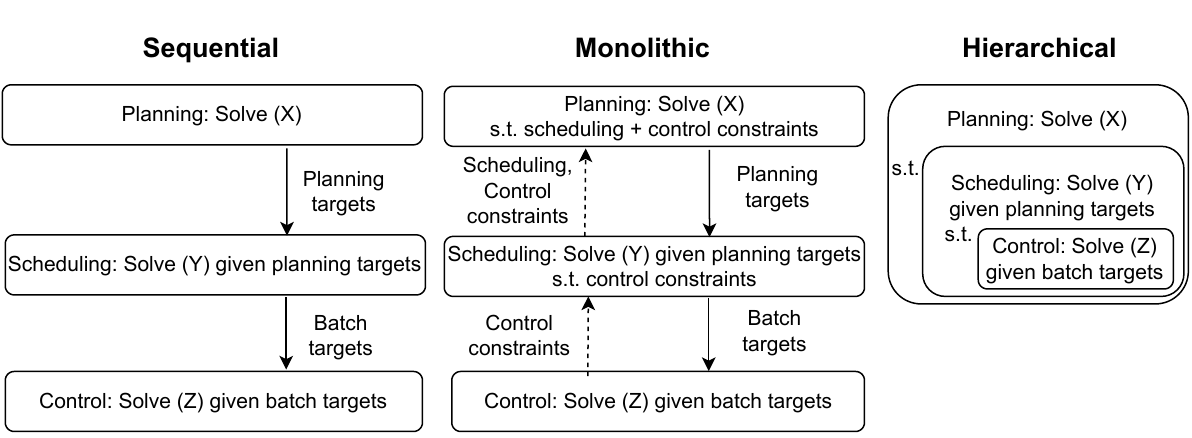}
    \caption{Hierarchical decision-making architectures: In the sequential approach, higher-level solutions feed into lower levels as setpoints. Each level solves their own optimization problem to meet these setpoints without explicit considerations of the lower levels. The monolithic approach resembles the sequential approach, but where feasibility in lower levels is satisfied via explicit incorporation of their constraints. The hierarchical approach could be viewed as a multi-level leader-follower game. A single `tri-level' optimization problem is solved, where the planning is subject to the optimal solution of the scheduling level, and the scheduling is constrained on the optimal solution of the control level.}
    \label{fig:int_frame}
\end{figure}

Companies within the process industries rely on mathematical optimization for their operations to remain competitive in an environment of increasingly stringent safety, environmental, and economic requirements \citep{Gounaris2019}. This gives rise to the field of enterprise-wide optimization (EWO) with the ultimate goal to coordinate all decision-making within a company \citep{Grossmann2005, Grossmann2012}. Significant value can be captured by integrating units across all hierarchical levels of decision-making (from design, planning, scheduling, to control). We adopt the terminology of \citet{CHU20152} and distinguish between the sequential, monolithic, and hierarchical integration of optimization models as shown in Figure \ref{fig:int_frame}. Conventionally, decision-making happens \textit{sequentially}: upper-level decisions are taken while disregarding lower-level considerations, and then fed as setpoints to the lower levels. There is however no guarantee that these setpoints are feasible in lower-level problems. Sequential decisions historically arise out of necessity due to tractability and organizational constraints.

In the \textit{monolithic} approach, lower-level feasibility considerations are included into upper-level optimization problems. This comes at the expense of a significant drop in computational tractability due to heterogeneity in model formulations and vastly different time horizons between levels. In the \textit{hierarchical} approach, upper-level decision-makers consider not only lower-level feasibility, but also optimality, accounting for how lower-level decisions affect upper-level objectives. This results in multi-level formulations, which are numerically intractable and mathematically difficult \citep{bihard}. 

There are two main tools that can be leveraged to alleviate the computational burden of integrated optimization problems. The first are relaxation or aggregation techniques \citep{newest} with the aim to relax the constraints, granularity, or detail of the integrated optimization models. This includes the replacement of some model parts - for example the lower level - by surrogate models \citep{TSAY201922, data, doi:10.1021/ie501986d} that are easier to handle by numerical solvers. Introducing surrogates into the original formulation inevitably incurs the risk of losing solution quality. As such, decision-makers inherently trade-off surrogate accuracy with optimization tractability when they choose the type of surrogate.
The second approach exploits the mathematical structure of integrated decision-making problems via decomposition algorithms \citep{decomp, MORAMARIANO2020106713}. Hierarchical problems could be viewed as coordination problems between one or multiple leader(s) and follower(s): each player has decision freedom over their respective problem but needs to coordinate on a sparse set of complicating variables. Traditional sequential decision architectures often lead to the coordinated variables to be few relative to the upper- or lower-level specific variables. As such, Lagrangean, and Benders decomposition \citep{KELLEY2018951, JI2021107166} for example lend themselves well to tractably solving hierarchically integrated problems.

The interest and successes in Machine Learning (ML) over the last decade have successfully diffused through the chemical engineering literature in recent years \citep{ML1,ML2}. Advances in ML have also inspired new developments in hierarchical planning, scheduling, and control. While surrogates - also known as meta-, scale-bridging, or reduced-order models \citep{10.1115/1.2429697, DU201559, TSAY201922} - have been established tools in the process systems engineering literature, ML has shifted attention to other approximation models \citep{BHOSEKAR2018250}: artificial neural networks, decision trees, Gaussian Processes, and support vector machines. We use `feasibility surrogates' and `optimality surrogates' as umbrella terms for any approximation model with the aim of mapping lower-level feasibility, or lower-level optimal response variables (objective or decision variables) of upper-level complicating variables or `setpoints'. This approach has analogues and similarities across many disciplines, from multi-parametric model predictive control and explicit optimal policies \citep{AVRAAMIDOU201917}, to amortized learning \citep{Amos2022TutorialOA} in ML and value functions \citep{Sachio} in Reinforcement Learning (RL).

The developments in ML have also heavily influenced the use of surrogates in derivative-free optimization (DFO) \citep{Larson2019}. DFO deals with the optimization of systems without having explicit or cheap access to gradient information. DFO algorithms are typically classified into either 'direct' methods that directly handle function evaluations, or 'model-based' methods that rely on the intermediate construction and optimization of surrogates. Although there are subtle differences, here, we use the term DFO synonymously with black-box, simulation-based, zeroth-order, or gradient-free optimization \citep{Amaran2016}. As DFO gains traction within the chemical engineering community \citep{vandeBerg2022Data-drivenApplications}, DFO has been exploited to solve integration problems that traditionally are targeted via decomposition algorithms. This includes coordination \citep{DDC} and multi-level problems in process systems engineering \citep{Zhao2021} and even in the hierarchical integration of process operations \citep{BEYKAL2022107551}. 

When the computational budget is available to attempt integration, we argue that data-driven techniques can be used to close the gap between monolithic and hierarchical approaches. In Section \ref{sec: lit rev}, we explore related works. We highlight how data-driven techniques fall into either the optimality surrogate or DFO approach and highlight novel aspects of our work. In Section \ref{sec: novelty}, we emphasize the aim of our work and novel aspects.

\subsection{Related works} \label{sec: lit rev}
There is ample review literature about solution strategies and applications in the general domain of integrated planning-scheduling-control (iPSC). Some works narrow in on the integration of scheduling and control (iSC) \citep{DAOUTIDIS2018179}, others focus on uncertainty \citep{DIAS201698}, their application to sequential batch processes \citep{7172222}, or the difference between top-down and bottom-up approaches \citep{BALDEA2014377, CASPARI202050}.
\citet{newest} present the most recent review at the time of writing and classify works according to their application to integrated planning and scheduling (iPS), iSC, or iPSC.  While most reviews are structured following the scope of integration, we are more interested in coming up with a typology of solution approaches. We discuss different solution methods and highlight their relationship to the general paradigms of surrogate modelling, decomposition, and derivative-free optimization (DFO). While we initially zone in on monolithic or simultaneous formulations, we argue that the same typology of solution approaches can be applied to hierarchical, multi-level formulations.

\subsubsection{Surrogates}
We use the term `surrogate models' to encompass all approximation models that reduce model complexity for the sake of tractability. Surrogates could be loosely regarded as relaxation or aggregation techniques and as such inevitably incur the risk of losing solution quality by forfeiting model accuracy. These approaches can be categorized into model reduction or system identification techniques \citep{TSAY201922}. While the former aims to derive a low-order model from detailed dynamic models (i.e. orthogonal decomposition etc.), system identification techniques learn a dynamic model from data. \citet{data} argue that data presents a natural bridge between different hierarchical layers and consequently advocate for the use of data-driven `time-scale bridging models' to integrate operations and control. Scale-bridging models are compatible with the `aggregate method' \citep{doi:10.1021/ie501986d}, where surrogates are used to map the linking functions between subproblems arising from the decomposition of different hierachical layers.

Surrogates differ by the amount of detail they capture from lower levels. As such, \citet{data} distinguish between static and dynamic models. Static surrogates for instance  do not explicitly account for transient data/dynamics. They could for instance learn only scheduling-relevant quantities from the control problem in iSC \citep{DU201559}. They can also be used to map the feasible space of lower levels \citep{XENOS2016418}. This is related to feasibility analysis \citep{BADEJO2022107759}. Convex region surrogates are especially popular as they map feasible sets using disjoint polytopes on historical operating data. Hence, the feasible region can be embedded into mixed-integer linear programming (MILP) formulations, avoiding the use of nonlinear terms and additional complexity. Nonlinear regression techniques could also be used as static models in mapping continuous outputs rather than binary feasibility: Gaussian Processes, support vector machines, piece-wise affine regression, or artificial neural networks would be well-suited to map the (optimal) cost of a lower level.

Dynamic surrogates explicitly account for time as an additional input. They are necessary when the dynamics of the lower level become relevant. Data-driven dynamic surrogates are overwhelmingly influenced by the system identification literature. Relevant techniques include Hammerstein-Wiener, linear state-space models \citep{KELLEY2022117468, TSAY20181273}, and Kalman filtering for online updates \citep{TSAY2021140}. 
The difference between learning dynamic surrogates or upper-level relevant variables divides the community: The system identification literature champions the learning of dynamic models from data which can subsequently be deployed within optimization-based feedback control such as model predictive control (MPC) \citep{BRADFORD2020106844}; The RL literature advocates for the approximation of optimal policies or feedback laws directly from data. \citet{spicyOne} argue that Q-Learning is not competitive with standard system identification on a simple, linear two air-zone heating, ventilation, and air conditioning system in the presence of small amounts of Gaussian noise. However, they find promise in the use of neural networks to learn the response of an MPC controller on a linear dynamic model with quadratic stage cost, and linear constraints.

\subsubsection{Decomposition}

Planning, scheduling, and control problems are often tackled by different business units within the same company, since they cover widely different timespans, types of decisions, constraints, and objectives. As such, the interactions between hierarchical decision-making layers, also called linking or complicating variables, naturally end up being sparse. Examples include planning targets and batch processing times in iPS or iSC problems respectively. When moving away from sequential to monolithic or even hierarchical integrated problem formulations, this gives rise to mathematical optimization structures that can be exploited by decomposition techniques \citep{decomp}. While these techniques describe a wide range of approaches, they typically consist of solving problems by breaking them up into smaller subproblems, which are solved separately in parallel or sequentially. This often involves many iterations over the subproblems coordinated by the solution of a `master problem'.

Lagrangean decomposition can be used to break up the iPSC problem into separate, more manageable iPS and control problems \citep{MORAMARIANO2020106713}. The coupling variables can even be learned via systematic network approaches. \citet{MITRAI202063} use community detection to find the mathematical structure in an iSC problem and the linking variables that can be exploited using Generalized Benders decomposition. As expected, they recover the scheduling problem in the first level and dynamic optimization subproblems in the second level.
In later work, \citet{MITRAI2022107859} use Multi-cut generalized Benders decomposition for iPSC by learning the problem structure and find that the planning-scheduling and control constraints are coupled only via the transition times. 

Decomposition can also be combined with surrogates. \citet{KELLEY2018951} use Hammerstein-Wiener models to integrate control dynamics into scheduling and solve the resulting formulation using Lagrangean decomposition. \citet{JI2021107166} use a combination of Generalized Benders Decomposition and genetic algorithms in the simultaneous solution of iSC.

\subsubsection{Derivative-free optimization}

Derivative-free optimization (DFO) is used to optimize systems without explicit gradient expressions \citep{Larson2019}. As such, DFO does not need access to the system's explicit model formulations. DFO leverages input-output evaluations, or `data', directly to find the optimum of the problem at hand through `direct' or model-based methods. Direct methods comprise of a wide range of techniques from random search, to grid search, and more sophisticated algorithms like the simplex method (Nelder-Mead\citep{Spendley1962}), adaptive mesh search (NOMAD \citep{NOMAD}) or DIviding REctangles (DIRECT \citep{DIRECT}) . Evolutionary methods \citep{Evol} could also be considered a sophisticated version of random and hence direct search. Developments in model-based DFO are closely related to surrogate techniques, as model-based DFO relies on the intermediate construction and optimization of approximation models to find the optimum. Popular methods range from more exploitative linear or quadratic trust region methods (COBYLA \citep{COBYLA} and BOBYQA \citep{BOBYQA}), to more explorative Bayesian Optimization using Gaussian Processes \citep{BayesOpt} or gradient-boosted trees \citep{Entmoot}.
DFO can exploit the same mathematical problem structure as decomposition approaches in iPSC. DFO approaches are especially prevalent in iPS, when the number of complicating variables (i.e. the planning targets) is few. Early attempts mostly use genetic algorithms \citep{9475992, 9226090} or surrogates/metamodels \citep{10.1115/1.2429697} under the name of simulation-based optimization \citep{WAN20041009}.

\subsubsection{DFO and single-level reformulation of multi-level problems}

 Accounting for game-theoretical considerations in the optimization of integrated supply chains and process systems often involves multi-level formulations \citep{10.3389/fceng.2023.1130568}. While sometimes simultaneous formulations for iPSC can be solved without previous solution approaches \citep{doi:10.1021/ie402563j}, the nested decision-making architecture of bi-level formulations require specialized approaches. \citet{Colson2007AnOO} discuss general bi-level optimization solution approaches.
Solving multi-level problems requires their reformulation into a form suitable for conventional single-level optimization. These techniques can be categorized into two main approaches: DFO and single-level reformulation. 
 
DFO extends to hierarchical and as such bi-level rather than simultaneous formulations by using data-driven optimization to find optimal planning or scheduling setpoints given the optimal response of the scheduling \citep{BEYKAL2022107551} or control level \citep{DIAS2018139} respectively. 
 On the other hand, single-level reformulation involves transforming the lower-level optimization problem into its Karush-Kuhn-Tucker (KKT) conditions. These can then be embedded into the upper level as big-M constraints, or using branching techniques on the complementarity conditions \citep{kleinert_schmidt_2023}. KKT reformulation has been used to embed linear MPC into scheduling \citep{SIMKOFF2019287}. This can also be applied to the integation of supply chain design and operation: all possible lower-level problems are enumerated, reformulated into their KKT conditions, embedded into the upper level, and finally solved using  decomposition \citep{YUE201781}. 
KKT reformulation is still computationally expensive even for linear cases and generally loses any theoretical guarantees in nonconvex iPSC.

\subsubsection{Surrogates for multi-level problems}

Rather than embed lower levels into the upper level using KKT reformulations, we can `learn' the optimal lower-level or follower response as a function of the leader's variables \citep{molan_schmidt_2023}. This could be considered an extension of the surrogate approach: rather than learn lower-level \textit{feasibility} from historical operating data, we can learn lower-level \textit{optimality} from the optimal solution of scheduling or control problems. We call `feasibility surrogates' and `optimality surrogates' the use of time-scale bridging models applied to simultaneous and hierarchical approaches respectively. The line between monolithic and hierarchical formulations can be blurry in feasibility analysis \citep{BADEJO2022107759} and often boils down to practical details - if feasible operating data is available or needs to be obtained via solving feasibility problems. 

The application of `optimality surrogates' in the integration of design, operations, and control appears in different domains under various names, from decomposition algorithms using `response functions' \citep{doi:10.1021/ie404272t}, to multi-parametric programming \citep{AvraamidouPistikopoulos, BURNAK2019164, PISTIKOPOULOS201685}. More recently, successes in Reinforcement Learning and Machine Learning have introduced new terminology and different perspectives. \citet{Sachio} have used Reinforcement Learning to learn an optimal policy on MPC simulations before embedding into the design problem. Software advances such as OMLT \citep{NEURIPS2021_17f98ddf, ceccon2022omlt, CECCON202257} are streamlining the integration of optimization and ML pipelines by automating reformulations of neural networks and decision trees and their subsequent embedding into Pyomo as constraint blocks that can be handled by MILP solvers. This could catalyze new advances in the safe deployment of learnt optimality surrogates to mitigate any model approximation risks, such as adversarial approaches or safe Reinforcement Learning \citep{SCHWEIDTMANN2019937, safeRLCalvin}.

\subsection{Aims and novelty} \label{sec: novelty}

While there is ample literature related to using optimality surrogates or DFO in solving integrated planning-scheduling-control problems, \citet{ESCAPE} are the first to attempt a tri-level solution of the iPSC problem by using DFO to optimize the planning level and an approximate scheduling-control surrogate. In this work, we strive to answer the questions that are left unanswered. We provide various methodologies that outperform the sequential approach on a small but realistic multi-site, multi-product integrated planning-scheduling-control case study:
\begin{itemize}
    \item We provide a tutorial-like exposition of the DFO and surrogate approaches and how we can evaluate solution qualities based on increasingly complex single- to tri-level evaluations.
    \item We investigate if we can integrate the scheduling-control layers as a surrogate into the planning and if this surrogate should be trained on scheduling-only, approximate, or exact scheduling-control data.
    \item We investigate if we can use DFO to find the tri-level solution without relying on model approximations.
    \item We explore practical techniques for balancing solution accuracy and computational time through combining the two methods.
\end{itemize}

In what follows, we illustrate our multi-site, multi-product, hierarchical planning-scheduling-control case study in Section \ref{sec: Case study}, and detail how the various levels are connected.  In Section \ref{sec: Methodology}, we present, in tutorial fashion, three combinations of the DFO and optimality surrogate approaches that can be leveraged towards a tractable tri-level planning-scheduling-control solution of our problem. We also highlight smart practical considerations that are crucial to making these approaches scalable. In Section \ref{sec: Results}, we present our results on the tractability-accuracy trade-off between the different methods on multiple solution quality and solution time metrics. We emphasize the role that understanding the interplay of ML and optimization pipelines has on navigating the inherent accuracy-tractability trade-off before concluding in Section \ref{sec: Conclusion}.

\section{Case study} \label{sec: Case study}

\begin{figure}[htp!] 
    \centering
    \includegraphics[width=\textwidth]{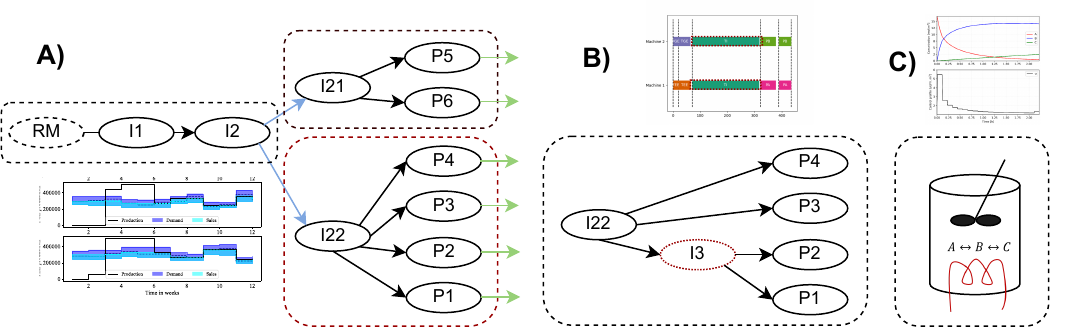}
    \caption{Outline of the hierarchical planning A), scheduling B), and control C) levels in the multi-site, multi-product case study. A) shows the planning state network denoting the materials ($\circ$) as well as processing ($\rightarrow$), transport ($\color{blue} \rightarrow \color{black}$), and sales ($\color{green} \rightarrow \color{black}$) across the three sites ($\square$). Scheduling and control are only relevant for the second site. B) shows the state network of the scheduling problem that includes an additional intermediate material. The production of I3 and P1-P4 needs to be scheduled in batches on 2 machines in 7 event time points. C) Each batch production is then subject to an optimal control problem.}
    \label{fig:Hierarchical}
\end{figure}

In this section, we present the planning, scheduling, and control layers from our multi-site, multi-product problem as shown in Figure \ref{fig:Hierarchical}. We then discuss how the decisions between the levels are interlinked as shown in Figure \ref{fig:iPSC}. We refer to Section \ref{sec: Appendix opt} for detailed optimization formulations of the high-level planning, scheduling, and control formulations presented in (\ref{eq: Planning}), (\ref{eq: Scheduling}), and (\ref{eq: Control}). The implementations of the planning, scheduling, and control formulations can also be found in the Github repository under \href{https://github.com/OptiMaL-PSE-Lab/DD-Hierarchical}{https://github.com/OptiMaL-PSE-Lab/DD-Hierarchical}.

\subsection{Planning}

The highest level of decision-making in our case study is planning. While planning model formulations can vary, production planning usually involves the solution of a (mixed-integer) linear program (MI)LP to determine the optimal production targets that optimize an economic objective over a specific planning horizon. Our planning problem (\ref{eq: Planning}) involves the solution of an LP over a planning horizon of 12 months. The state network of the planning is illustrated in Figure \ref{fig:Hierarchical}A. We consider three sites, where in the first site, raw material RM (assumed unlimited) undergoes two intermediate processing stages (I1 and I2). I2 is then transported to the two other sites, denoted as I21 and I22. In the second and third site, intermediate material I22 and I21 undergo further processing into one of four or one of two products (P1-P4 and P5-P6) respectively, which are then sold to meet external customer demand. The goal of the planning is to determine the production of each material given processing yields (black arrows), the material transport given lead times (blue arrows), and the sales given customer demand (green arrows) that maximize sales and minimize transportation, storage, and production cost. The LP is subject to inventory material balances on each material, site-wide resource utilization limits, and upper and lower bounds on inventory for each material. The full formulation can be found in (\ref{eq: app Planning}).

\begin{equation}
\label{eq: Planning}
\begin{aligned}
     & \min.  
     && \text{ Transport + Storage + Production - Sales } \\
     & \text{s.t.} 
     &&  \text{Initial conditions}\\
     &&& \text{Inventory mass balances with production and transportation}  \\
     &&& \text{Resource limits} \\
     &&& \text{Sales limits} \\
     &&& \text{Safe storage constraints} \\
\end{aligned}
\end{equation}

\subsection{Scheduling}

The optimal planning targets for each of the 12 planning timesteps are then fed to the scheduling layer. At each planning timestep, a scheduling problem determines the resource allocation required to reach the planning targets in minimal time, i.e. for each of the 7 scheduling timesteps (events), for each of the two machines, which production (job) happens and for how long. Our scheduling formulation is based on \citet{scheduling} and relies on the state-task network approach with a common continuous-time representation for all units. The formulation also accounts for variable batch sizes, variable processing times and sequence-dependent changeover times. In our integrated framework, the scheduling layer involves the solution of a MILP (\ref{eq: Scheduling}) at each planning step involving discrete variables in the assignment constraints (if a batch production happens in a specific machine at a given time) and changeover constraints (if we switch from one production to another at a given event time point and machine) on top of the planning-level mass balances. The batch constraints, production recipes, production duration, and changeover duration are machine- and sequence-specific. There is also a more accurate estimate of the inventory and resource limits available compared to the planning layer. The scheduling layer is only considered for site 2 given the simplicity of the state network in the other sites. Explicitly accounting for scheduling in the other sites is compatible with our proposed methodologies at a potential increase in solution time, especially since the scheduling problems are only loosely connected via few planning variables. The state-task network corresponding to site 2 is depicted in Figure \ref{fig:Hierarchical}B, where on top of the state network from Figure \ref{fig:Hierarchical}A, we also consider a perishable, bottleneck intermediate I3 between I22 and products P1 and P2. I3 is produced on-demand for each planning timestep with any unused production after the planning period going to waste. The full formulations can be found in (\ref{eq: app Scheduling obj}) to (\ref{eq: app Scheduling tightening}).

\begin{equation}
\label{eq: Scheduling}
\begin{aligned}
     & \min.  
     && \text{ Makespan } \\
     & \text{s.t.} 
     && \text{Inventory mass balances with production}  \\
     &&& \text{Batch constraints} \\
     &&& \text{Duration constraints} \\
     &&& \text{Assignment constraints} \\
     &&& \text{Changeover constraints} \\
     &&& \text{Resource limits} \\
     &&& \text{Meet monthly planning production targets} \\     
\end{aligned}
\end{equation}

\subsection{Control}

The control formulations are based on \citet{control}. Each batch control takes the batch targets for each job-machine allocation as determined in the scheduling layer, and determines the optimal cooling flowrate that achieves said batch target while minimizing a combination of processing time and energy cost. This involves the solution of an optimal control problem (\ref{eq: Control}) subject to final time quality conditions and nonlinear differential expressions of the batch kinetics and energy transfer. The full formulation can be found in (\ref{eq: app Control}).

\begin{equation}
\label{eq: Control}
\begin{aligned}
     & \min.  
     && \text{Processing time \& energy cost} \\
     & \text{s.t.} 
     && \text{Mass balances given by differential equations}  \\
     &&& \text{Kinetics dependent on state and control variables} \\
     &&& \text{Energy cost dependent on control variables} \\
     &&& \text{Final time quality conditions} \\
\end{aligned}
\end{equation}

\begin{figure}[htp!] 
    \centering
    \includegraphics[width=\textwidth]{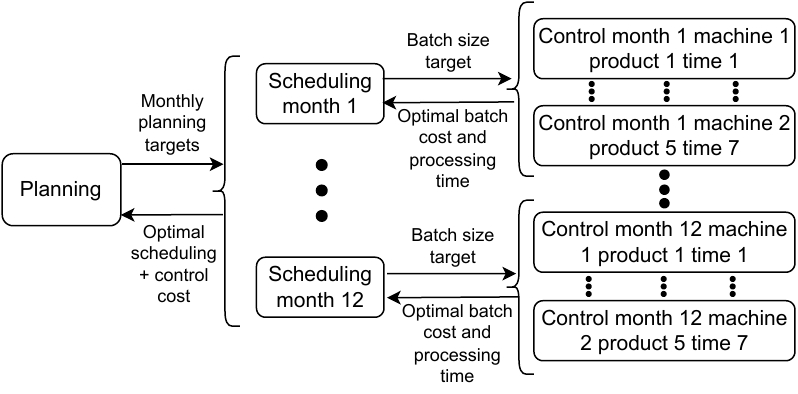}
    \caption{Linking variables between the planning, scheduling, and control levels. Each planning problem feeds the planning targets for its 12 timesteps to separate scheduling problems to be solved in parallel. Each of the 12 scheduling problems, in deciding which of the 5 productions to schedule on the 2 machine at the 7 events, obtains the optimal batch processing time and cost as a function of the batch size targets for each of the 70 possible batch assignments. The optimal scheduling and control finally feed back into the planning.}
    \label{fig:iPSC}
\end{figure}

\subsection{Integration} \label{sec: Integration}

iPSC problems are either solved sequentially by ignoring all lower-level considerations, monolithically by only accounting for lower-level feasibility, or hierarchically by explicitly accounting for lower-level optimality as depicted in Figure \ref{fig:int_frame}. While the sequential and monolithic approaches are often employed for tractability or practical reasons, the hierarchical approach most accurately reflects organizational decision-making. Figure \ref{fig:iPSC} shows the information flow between the layers. These connecting or complicating variables link the three decision-making units that can otherwise be solved autonomously by each level. At the highest level, the 12 optimal monthly planning variables feed as setpoints into a separate scheduling problem. For each planning timestep, all batch assignment targets as determined in the scheduling layer then feed as setpoints into the control level. In the scheduling layer, we need to decide on and as such explicitly consider all 5 possible productions for each of the 14 (7 events by 2 machines) optimal control problems. If we want to evaluate the scheduling layer for a given batch assignment however, we only need to solve for the 14  optimal control problems with the assigned production. Each optimal control solution then feeds back to the optimal processing times and corrects the makespan corresponding to a given sequence. The optimal changeover and cooling costs from the scheduling and control layer then feed back into the planning level where they correct the economic objective.

While it is in principle possible to solve a monolithic formulation of the iPSC, this is intractable in our case without employing surrogates or decomposition. First, the nonlinear dynamics corresponding to up to 70 job-machine assignments would have to be embedded into the scheduling, exploding the dimensionality of the scheduling and turning the integrated problem mixed-integer \textit{nonlinear}. Then, the changeover, assignment, and duration constraints of these 12 integrated problems would have to be embedded into the planning. 

A hierarchical solution is not just intractable but mathematically difficult given the nested layers of decision-making. Each tri-level optimization requires the solution of up to 840 control problems: the top layer (planning) is subject to the optimal solution of 12 scheduling optimization problems that themselves are constrained on the optimal solution of up to 70 optimal control problems each.

Data-driven techniques have become essential to the solution of iPSC formulations. In the next section, we want to show how we can exploit data-driven techniques to solve hierarchical (multi-level) formulations. While it is unrealistic to expect to find the optimal solution to tri-level optimization formulations (even the simplest linear-linear bilevel optimization problems are NP-hard \citep{bihard}), our aim is to find a tractable solution to the tri-level hierarchical formulation that outperforms the conventional sequential solution. In this work, we prioritize finding an approximate solution to an accurate model rather than finding the optimal solution to an approximate model. We acknowledge however that there are many reasons to use the sequential or monolithic approach, i.e. when the lower-level objectives are aligned with upper-level objectives, or when any potential gain in the tri-level formulation does not justify the increase in compute. 

\section{Methodology} \label{sec: Methodology}

To streamline the exposition of how we can use surrogates and derivative-free optimization to solve our iPSC problem, we first show how we can use both techniques in solving general bi-level problems.

\subsection{Using data to solve a canonical bi-level problem} \label{eq: canonical}

Problem (\ref{eq: Canonical}) presents a general bi-level formulation where the leader determines the set of variables that appear only in their level $\mathbf{x}_{up}$, and the connecting variables $\mathbf{x}_{u \rightarrow l}$ that appear as setpoints in the lower level. Given the linking variables $\mathbf{x}_{u \rightarrow l}$, the follower optimizes their objective $f_{low}(\cdot)$ by manipulating the variables that are specific to the lower level $\mathbf{x}_{low}$ and the connecting variables that also appear in the upper level $\mathbf{x}_{l \rightarrow u}$. As such, the leader optimizes their objective $f_{up}(\cdot)$ while explicitly accounting for the optimal response of the follower $\mathbf{x}^*_{l \rightarrow u}$ to $\mathbf{x}_{u \rightarrow l}$. We consider the complicating variables $\mathbf{x}_{u \rightarrow l}$ and $\mathbf{x}_{l\rightarrow u}$ to be continuous, but not necessarily $\mathbf{x}_{up}$ and $\mathbf{x}_{low}$. No such restrictions need to be placed on any of the variables but continuous variables are generally easier to handle.

\begin{subequations}
\label{eq: Canonical}
    \begin{align}
            & \underset{\substack{\mathbf{x}_{up}, \mathbf{x}_{u \rightarrow l}, \\ \mathbf{x}^*_{low}, \mathbf{x}_{l \rightarrow u}^*}}{ \text{min.}}
            && f_{up}(\mathbf{x}_{up}, \mathbf{x}_{u \rightarrow l},  \mathbf{x}_{l \rightarrow u}^* )\\
            & \text{s.t.}
            && \mathbf{h}_{up}(\mathbf{x}_{up}, \mathbf{x}_{u \rightarrow l},  \mathbf{x}_{l \rightarrow u}^* ) = \mathbf{0} \\
            &&& \mathbf{g}_{up}(\mathbf{x}_{up}, \mathbf{x}_{u \rightarrow l},  \mathbf{x}_{l \rightarrow u}^* ) \leq \mathbf{0} \\
            &&&  \mathbf{x}^*_{low}, \mathbf{x}_{l \rightarrow u}^* \in \quad \begin{aligned}[t]
            &  \underset{\mathbf{x}_{low}, \mathbf{x}_{l \rightarrow u}}{ \arg \min.} 
            && f_{low}(\mathbf{x}_{low}, \mathbf{x}_{l \rightarrow u},  \mathbf{x}_{u \rightarrow l} ) \\
            & \text{s.t.} 
            && \mathbf{h}_{low}(\mathbf{x}_{low}, \mathbf{x}_{l \rightarrow u},  \mathbf{x}_{u \rightarrow l} ) = \mathbf{0} \\
            &&& \mathbf{g}_{low}(\mathbf{x}_{low}, \mathbf{x}_{l \rightarrow u},  \mathbf{x}_{u \rightarrow l} ) \leq \mathbf{0} \\
            \end{aligned} \label{eq: low}
    \end{align}
\end{subequations}

We now present two different ways that data-driven techniques can be exploited towards a solution of bi-level problems.

\subsubsection{Optimality surrogates}

We first present optimality surrogates as a solution approach to multi-level problems. Crucial to solving bi-level optimization problems is their reformulation into a single level, such that the problem can be exploited using standard optimization solvers. In the optimality surrogate approach, we first sample various combinations of the complicating variables $\mathbf{x}_{u \rightarrow l}$ and solve the corresponding lower-level optimization problem to extract the optimal response $\mathbf{x}^*_{l \rightarrow u}$. Then, we can use any supervised learning technique to construct a model - neural networks $\mathcal{NN}( \cdot)$ in our case - to map  $\mathbf{x}^*_{l \rightarrow u}$ to $\mathbf{x}_{u \rightarrow l}$. The explicit expression $\mathbf{x}_{l \rightarrow u}^* = \mathcal{NN}( \mathbf{x}_{u \rightarrow l})$ can then be embedded as a constraint into the upper level, replacing the lower-level optimization problem, and collapsing the bi-level into a single-level formulation (\ref{eq: optimality}) that can be readily implemented in standard optimization software. 

\begin{equation}
    \label{eq: optimality}
    \begin{aligned}
            & \underset{\substack{\mathbf{x}_{up}, \mathbf{x}_{u \rightarrow l}, \\  \mathbf{x}_{l \rightarrow u}^*}}{ \text{min.}}
            && f_{up}(\mathbf{x}_{up}, \mathbf{x}_{u \rightarrow l},  \mathbf{x}_{l \rightarrow u}^* )\\
            & \text{s.t.}
            && \mathbf{h}_{up}(\mathbf{x}_{up}, \mathbf{x}_{u \rightarrow l},  \mathbf{x}_{l \rightarrow u}^* ) = \mathbf{0} \\
            &&& \mathbf{g}_{up}(\mathbf{x}_{up}, \mathbf{x}_{u \rightarrow l},  \mathbf{x}_{l \rightarrow u}^* ) \leq \mathbf{0} \\
            &&& \mathbf{x}_{l \rightarrow u}^* = \mathcal{NN}( \mathbf{x}_{u \rightarrow l} ) \text{\quad\quad\quad (neural network)}
    \end{aligned}
\end{equation}

In principal, it is possible to construct surrogates that present accurate `explicit control laws' for certain kinds of parameterized optimization problems \citep{rawlings2017model}. However, parametric programming explodes in complexity with the size of the lower levels \citep{PAROC}. Consequently, we train approximate optimality surrogates, forfeiting theoretical guarantees that the lower-level is solved to optimality. As such, care should be taken in choosing the type and architecture of the surrogate when trading off solution accuracy and tractability.

We suggest using neural networks with piecewise linear activation functions (ReLU) as optimality surrogates since their expressions can be reformulated into mixed-integer \textit{linear} constraints \citep{huchette2023deep}, avoiding the introduction of nonlinear terms into upper level objectives. Other techniques like decision trees could be used to the same end. There is ongoing debate as to whether it is favourable to use mixed-integer linear formulations with discrete solvers or to use global solvers on the full-space formulations of the surrogates \citep{SCHWEIDTMANN2019937}.

We note that the concept of optimality (or feasibility) surrogates can also be used in monolithic formulations. For example, 
if (binary) data on the feasible operating region is available, we can map the constraints which can be integrated into the upper level as $\mathcal{NN}( \mathbf{x}_{u \rightarrow l}) \leq 0$.

\subsubsection{Derivative-free optimization}

In our second approach, we leverage recent advances in derivative-free optimization (DFO). The idea of using derivative-free optimization (DFO) for bi-level formulations is to fix as many upper-level variables as necessary such that the upper level is fully determined. First, we observe that when $ \mathbf{x}_{u \rightarrow l}$ is fixed, we can obtain $ \mathbf{x}^*_{l \rightarrow u}$ by solving the lower-level optimization problem. On top of this, we can split up the upper level variables $\mathbf{x}_{up}$ into a set of decision variables 
$\mathbf{x}_{DFO}$ and $\mathbf{x}_{sim}$, such that $\mathbf{x}_{sim}$ can be obtained by solving the system of equality constraints $\mathbf{h}_{up}(\cdot) = \mathbf{0}$ at $\mathbf{x}_{DFO}$ and $\mathbf{x}_{u \rightarrow l}$ fixed. As such, we can use DFO to find $\mathbf{x}_{u \rightarrow l}$ and $\mathbf{x}_{DFO}$ that minimize the black-box objective consisting of the upper-level objective $f_{up}(\cdot )$ augmented by a penalization of any violation of $\mathbf{g}_{up}(\cdot) \leq \mathbf{0}$. The extent of penalization can be tuned via the penalty parameter $\rho$, which as a rule of thumb can be chosen to be an order of magnitude higher than the objective terms $f_{up}(\cdot)$. We essentially solve the general bi-level formulation (\ref{eq: Canonical}) as a single-level DFO problem (\ref{eq: DFO}) where each evaluation consists in the expensive solution of the lower-level optimization problem $\mathcal{LOW}(\cdot)$ equivalent to (\ref{eq: low}).

\begin{equation}
    \label{eq: DFO}
    \begin{aligned}
            & \underset{ \mathbf{x}_{DFO} , \mathbf{x}_{u \rightarrow l}}{ \text{min.}}
            && \mathcal{BB}(\mathbf{x}_{DFO}, \mathbf{x}_{u \rightarrow l}) \text{\quad\quad\quad (Black-box)} \\
            & \text{where}
            && \mathbf{x}_{l \rightarrow u}^* \leftarrow \mathcal{LOW}( \mathbf{x}_{u \rightarrow l} ) \text{\quad\quad (lower problem)} \\
            &&&  \mathbf{x}_{sim}   \leftarrow \mathbf{h}_{up}(\mathbf{x}_{sim}, \mathbf{x}_{DFO}, \mathbf{x}_{u \rightarrow l},  \mathbf{x}_{l \rightarrow u}^* ) = \mathbf{0} \\
            &&&  \mathbf{x}_{up} = [\mathbf{x}_{sim}, \mathbf{x}_{DFO}]   \\
            &&&  \text{penalty} = \rho || \max(\mathbf{0} , \mathbf{g}_{up}(\mathbf{x}_{up}, \mathbf{x}_{u \rightarrow l},  \mathbf{x}_{l \rightarrow u}^* ))||^2 \\
            &&& \mathcal{BB}(\mathbf{x}_{DFO}, \mathbf{x}_{u \rightarrow l}) = f_{up}(\mathbf{x}_{up}, \mathbf{x}_{u \rightarrow l},  \mathbf{x}_{l \rightarrow u}^* ) +  \text{penalty}\\
    \end{aligned}
\end{equation}

The left arrows $\leftarrow$ denote that we obtain $\mathbf{x}_{sim}$ and $ \mathbf{x}_{l \rightarrow u}^*$ from the solution of the set of equality constraints $\mathbf{h}_{up}(\cdot) = \mathbf{0}$ and of the lower-level optimization instance $\mathcal{LOW}(\cdot)$ within the black-box simulation respectively. In principle, we can use any kind of DFO solver for the single-level DFO problem (\ref{eq: DFO}). However, since each evaluation involves the solution of an optimization problem, this severely restricts the available evaluation budget especially since these problems tend to be higher-dimensional than the applications at which DFO excels. This prevents the practical applicability of many over-explorative methods like evolutionary search, and Bayesian Optimization. We suggest using exploitative trust-region methods and leveraging as much problem knowledge as possible in finding a good initial guess.

\subsubsection{Differences between the two methods}
Model-based DFO and the surrogate approach could be easily confused as they both rely on surrogates. The surrogate approach constructs relevant input-output models of the lower level before embedding these surrogates into the upper level, where a single optimization formulation is solved using algebraic modelling languages such as Pyomo or GAMS. In model-based DFO, surrogates are only used in the `outer level' to trade-off exploitation and exploration and determine the next sample input to the simulation. Although we have opted for model-based DFO methods in this work, direct DFO methods could be used (e.g., simplex, metaheuristics) instead requiring no surrogates whatsoever. In the `inner level', the simulation then calls the original lower-level optimization problem (e.g. in Pyomo \citep{pyomo}) as a black-box to return relevant quantities to the upper level.

More importantly, both approaches present vastly different model accuracy versus solution tractability trade-offs. While the surrogate approach inevitably runs the risk of losing solution quality through model inaccuracies, DFO scales poorly with the expense of lower-level problems and the number of DFO dimensions. This motivates a careful investigation into the intricacies of our case study and how these can be used to integrate the planning-scheduling and the scheduling-control layer.

\subsection{Hierarchical planning-scheduling-control as tri-level formulation}
Let us formally present the mathematical problem we are addressing. 
Problem (\ref{eq: tri-level}) formalizes the hierarchical iPSC formulation described in Section \ref{sec: Integration}. The upper level consists of LP formulations to determine the planning-specific and complicating planning target variables ($\mathbf{x}_p$ and $\mathbf{x}_{prod}$) that minimize the planning-level economic objective $f_p(\cdot)$ augmented by the optimal scheduling and control costs $c^*_s$ and $c^*_c$. $c^*_s$ is obtained by constraining the upper level on solving 12 MILP problems. Their aim is to find the scheduling-specific and batch target variables ($\mathbf{x}_s$ and $\mathbf{x}_{batch}$) that minimize the makespan required to fulfil $\mathbf{x}_{prod}$. The optimal batch processing times $ t^*_f$ however rely on optimal control operation. Each of the 12 scheduling problems are subject to ($14-70$) nonlinear optimal control problems that determine the control-specific variables $\mathbf{x}_c$ and processing times $t_f$ to meet the batch targets ($\mathbf{x}_{batch}$) that minimize a mixture of energy cost and processing time $ f_c(\cdot)$.

\begin{equation}
    \label{eq: tri-level}
    \begin{aligned}
            & \underset{\substack{\mathbf{x}_p, \mathbf{x}_{prod}, \\  \mathbf{x}^*_s, c_s^*,  \mathbf{x}_b^*,  \\  c^*_c,  t^*, \mathbf{x}^*_c}}{ \text{min.}}
            && f_p(\mathbf{x}_p, \mathbf{x}_{prod}) + c^*_s +  c^*_c  \\
            & \text{s.t.}
            && \mathbf{h}_p(\mathbf{x}_p, \mathbf{x}_{prod}) = \mathbf{0} \\
            &&& \mathbf{g}_p(\mathbf{x}_p, \mathbf{x}_{prod}) \leq \mathbf{0} \\
            &&&  \substack{ \mathbf{x}^*_s, \mathbf{x}_b^*, c_s^*, \\  \mathbf{x}^*_c, t^*, c^*_c} \in \quad \begin{aligned}[t]
            &  \underset{\substack{ \mathbf{x}_s,  \mathbf{x}_{batch},  c_s, \\  \mathbf{x}^*_{ctrl}, t_f^*, c^*_{ctrl}}}{ \arg \min.} 
            && f_s(\mathbf{x}_s, \mathbf{x}_{batch},  t_f^*,  \mathbf{x}_{prod} ) \\
            & \text{s.t.} 
            && \mathbf{h}_s(\mathbf{x}_s, \mathbf{x}_{batch}, c_s,  t_f^*,  \mathbf{x}_{prod} ) = \mathbf{0} \\
            &&& \mathbf{g}_s(\mathbf{x}_s, \mathbf{x}_{batch},  t_f^*,  \mathbf{x}_{prod} ) \leq \mathbf{0} \\
            &&&  x^*_{ctrl}, c^*_{ctrl}, t^*_f \in \quad \begin{aligned}[t]
            & \underset{\mathbf{x}_c, c_{ctrl}, t_f}{ \arg \min.} 
            && f_c(\mathbf{x}_c, t_f ,  \mathbf{x}_{batch} ) \\
            & \text{s.t.} 
            && \mathbf{h}_c(\mathbf{x}_c,  t_f, c_{ctrl},  \mathbf{x}_{batch} ) = \mathbf{0} \\
            &&& \mathbf{g}_c(\mathbf{x}_c, t_f,  \mathbf{x}_{batch} ) \leq \mathbf{0} \\
            \end{aligned} 
            \end{aligned} 
    \end{aligned}
\end{equation}

Before we describe how we can apply the techniques described as applied to the general bi-level formulation (\ref{eq: canonical}), we first determine four metrics by which we can assess the quality of a planning-level solution in increasing levels of complexity and accuracy.

\subsection{Metrics}

\begin{figure}[htp!] 
    \centering
    \includegraphics[width=\textwidth]{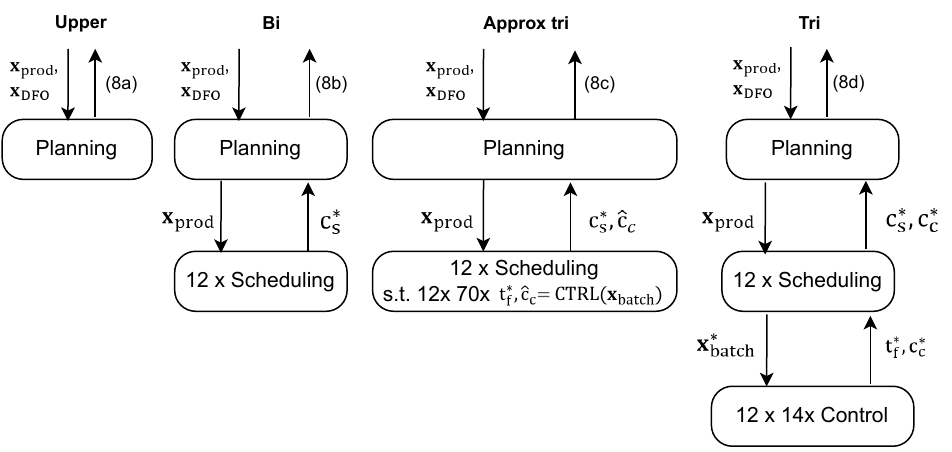}
    \caption{Increasing levels of solution accuracy evaluation. The upper evaluation only evaluates the planning objective and constraints in a simulation. The bi evaluation solves 12 scheduling problems in parallel in each evaluation to account for the scheduling cost in the objective. The approximate tri evaluation embeds control surrogates into the scheduling to account for approximate control cost in the objective. The tri evaluation solves scheduling and control levels in sequence to emulate hierarchical decision-making and adds the `real' hierarchical scheduling and control costs to the planning objective.}
    \label{fig:Evaluations}
\end{figure}

To compare the solution quality of the proposed solution approaches, we define four different solution metrics. Figure \ref{fig:Evaluations} illustrates the difference between the four solution quality evaluation metrics accounting for different levels of integration of lower-level information into the upper-level planning objective (\ref{eq: metrics}). In each case, we start by extracting the $\mathbf{x}_{DFO}$ and $\mathbf{x}_{prod}$ variables, used as input to the DFO (Section \ref{sec: meth1}) or obtained from the solution of the conventional optimization instance in the upper level (Section \ref{sec: meth2}). We explain in the next section how the planning-level variables are partitioned into $\mathbf{x}_{DFO}$. This information is sufficient in simulating the planning level and by extension hierarchical solutions by feeding down optimal solutions sequentially. These evaluations constitute a crucial part in applying the DFO method (\ref{eq: DFO}) on the planning-scheduling integration and is presented in greater detail in the next section.
 
 The upper or \textit{planning-only} evaluation (\ref{eq: metrics1}) only simulates the planning-level objective and any upper-level inequality constraint penalizations $f_p(\cdot) + g_{viol}$. The other metrics differ in how much detail is captured in the optimization at each evaluation: The \textit{bi} evaluation (\ref{eq: metrics2}) calls all scheduling optimization problems in parallel to augment the upper objective by $c^*_s$; \textit{approx tri} (\ref{eq: metrics3}) calls integrated scheduling-control optimization formulations instead to obtain the optimal scheduling and approximate optimal control costs $c^*_s + \hat{c}_c$; \textit{tri} (\ref{eq: metrics4}) solves all scheduling problems first whose optimal solutions feed into the optimal control to obtain the `real' tri-level objective using $c^*_s + c^*_c$. These evaluation metrics present increasingly accurate hierarchical decision-making modelling with the upper (\ref{eq: metrics1}) and tri (\ref{eq: metrics4}) evaluations emulating the planning decisions according to the sequential and hierarchical frameworks respectively as depicted in Figure \ref{fig:int_frame}.

\begin{subequations}
\label{eq: metrics}
    \begin{align}
            & \text{Planning-only:}
            && f_{up}(\cdot) + g_{viol} \label{eq: metrics1}\\
            & \text{Bi:}
            && f_{up}(\cdot) + g_{viol} + c^*_s \label{eq: metrics2}\\
            & \text{Approx tri:}  
            && f_{up}(\cdot) + g_{viol} + c^*_s + \hat{c}_c \label{eq: metrics3}\\
            & \text{Tri:} 
            && f_{up}(\cdot) + g_{viol} + c^*_s + c^*_c \label{eq: metrics4}
    \end{align}
\end{subequations}

\begin{figure}[htp!] 
    \centering
    \includegraphics[width=\textwidth]{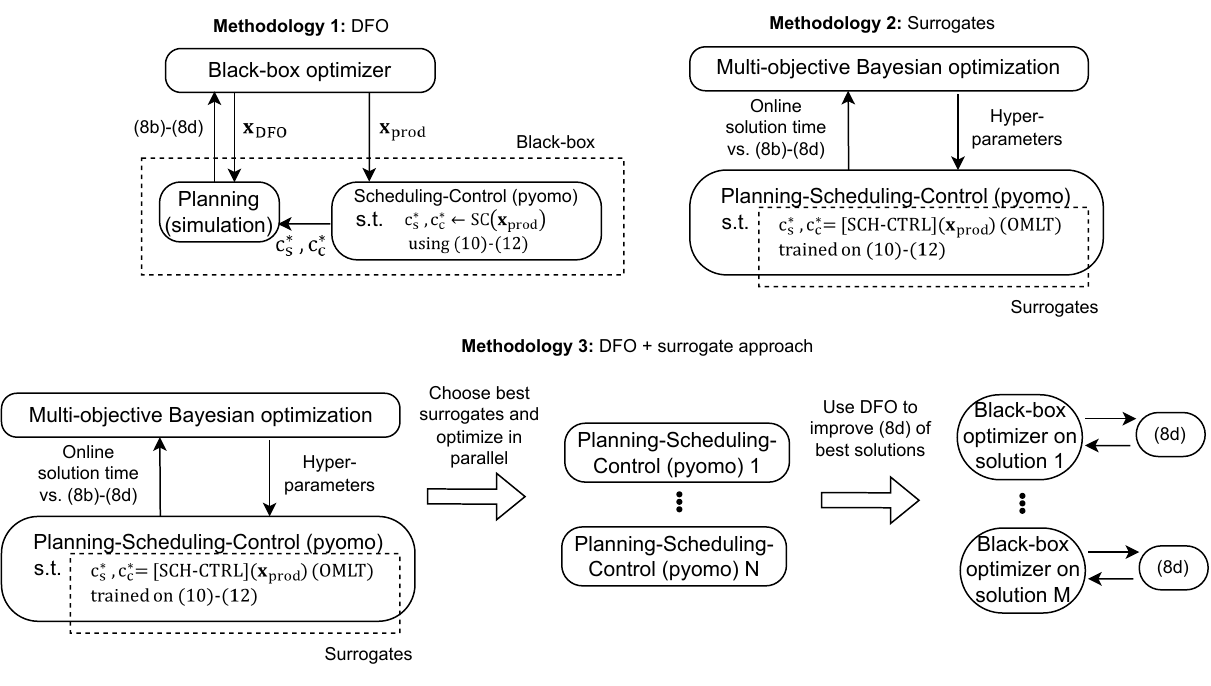}
    \caption{The three proposed methodologies: The first methodology uses DFO to optimize any of the solution metrics in (8b-8d) directly by considering the planning simulation and scheduling-control optimization as a black-box. The second methodology uses hyperparameter optimization to find a surrogate architecture that, trained on one of the solution metrics in (8b-8d), leads to the best solution accuracy-time trade-off after embedding into the planning. The third methodology combines both workflows in using DFO to optimize the solutions obtained from the surrogate approach.}
    \label{fig:Methods}
\end{figure}

 These metrics become crucial in all three methodologies as described in Figure \ref{fig:Methods}: the DFO, surrogate, and combined approach. In principle, any solution quality metrics (\ref{eq: metrics2})-(\ref{eq: metrics4}) can be used as the DFO objective; Similarly, we can map surrogates to map the outputs of scheduling-only (\ref{eq: bi eval}), approximate scheduling-control (\ref{eq: approx tri eval}), or accurate scheduling-control (\ref{eq: tri eval}) optimization problems called in the bi, approx tri, and tri evaluation metrics (\ref{eq: metrics2})-(\ref{eq: metrics4}). A clear distinction between the increasing levels of solution evaluation becomes important in presenting the results. In the next section, we illustrate how we use DFO to find the planning-level variables that optimize any one of the solution metrics from (\ref{eq: metrics}).

\subsection{DFO in the planning level} \label{sec: DFO init}

In order to use DFO for the integration of planning and scheduling, we reformulate the tri-level problem (\ref{eq: tri-level}) into the canonical DFO formulation (\ref{eq: DFO}) for bi-level problems. We split the planning-specific variables $\mathbf{x}_{p}$ of the tri-level optimization problem (\ref{eq: tri-level}) into its constituent storage $\mathbf{x}_{store}$, transport $\mathbf{x}_{transp}$, and sales $\mathbf{x}_{sales}$ variables.  We observe that given $\mathbf{x}_{prod}$, $\mathbf{x}_{transp}$, and $\mathbf{x}_{sales}$, we can exploit the inventory mass balances $\mathbf{h}_{p}(\cdot)$ to obtain $\mathbf{x}_{store}$. We now have all the variables required to compute the quantities of interest that appear in the upper level: $f_p(\cdot)$, the penalty term $g_{viol}$ via $\mathbf{g}_p(\cdot)$, and any (approximate) optimal costs $c^*_s$ and $c^*_c$ from simulating lower-level optimization problems via $\mathcal{SC}(\cdot)$. Problem (\ref{eq: DFO only}) shows the formulations involved in using DFO to find $\mathbf{x}_{sales}, \mathbf{x}_{transp}, \mathbf{x}_{prod}$ that optimize the integrated planning-scheduling-control problem as a black-box objective $\mathcal{BB}(\cdot)$.

\begin{equation}
    \label{eq: DFO only}
    \begin{aligned}
            & \underset{\mathbf{x}_{sales}, \mathbf{x}_{transp}, \mathbf{x}_{prod}}{ \text{min.}}
            && \mathcal{BB}(\mathbf{x}_{sales}, \mathbf{x}_{transp}, \mathbf{x}_{prod}) \\
            & \text{where}
            && \mathbf{x}_{store} \leftarrow \mathbf{h}_p(\mathbf{x}_{store}, \mathbf{x}_{sales}, \mathbf{x}_{transp}, \mathbf{x}_{prod}) = \mathbf{0} \\
            &&& \mathcal{BB}(\cdot) \leftarrow f_p(\mathbf{x}_{store}, \mathbf{x}_{sales}, \mathbf{x}_{transp}, \mathbf{x}_{prod}) + g_{viol} +  c^*_s +  c^*_c   \\
            &&& g_{viol} \leftarrow \mathbf{g}_p(\mathbf{x}_{store}, \mathbf{x}_{sales}, \mathbf{x}_{transp}, \mathbf{x}_{prod}) \leq \mathbf{0} \\
            &&&  c_s^*,  c^*_c  \leftarrow \mathcal{SC}( \mathbf{x}_{prod} ) \\ 
    \end{aligned}
\end{equation}

Comparing (\ref{eq: DFO only}) with the canonical DFO formulation (\ref{eq: DFO}), we see that $\mathbf{x}_{DFO} = [\mathbf{x}_{transp}, \mathbf{x}_{sales}]$, $\mathbf{x}_{sim}=\mathbf{x}_{store}$, $\mathbf{x}_{u \rightarrow l} = \mathbf{x}_{prod}$, and $\mathcal{LOW}(\cdot) = \mathcal{SC}(\cdot)$. In practice, this means that the DFO solver determines only the complicating planning targets, as well as a subset of the planning-specific variables (the transport and sales variables). This information is sufficient in determining all other relevant quantities in computing the black-box planning-scheduling-control objective. Different formulations of the (integrated) scheduling-control optimization problem $\mathcal{SC}(\cdot)$ lead to different solution accuracy metrics (\ref{eq: metrics}) such that $\mathbf{x}_{l \rightarrow u} = \emptyset, c^*_s, c^*_s+\hat{c}_c, \text{ or } c^*_s+c^*_c,$. \textit{Planning-only} does not use $\mathcal{SC}(\cdot)$. In the \textit{bi} metric evaluations, $\mathcal{SC}(\cdot)$ solves the 12 scheduling optimization problems taking the form of (\ref{eq: bi eval}) to extract $c^*_s$ without any restrictions on $t_f$.  

\begin{equation}
\label{eq: bi eval}
    \begin{aligned}
        & \mathcal{SC}( \mathbf{x}_{prod} ) =  &&  \underset{ \mathbf{x}_s,  \mathbf{x}_{batch},  c_s}{ \arg \min.} 
            && f_s(\mathbf{x}_s, \mathbf{x}_{batch},  t_f,  \mathbf{x}_{prod} ) \\
            & \text{ } && \text{s.t.} 
            && \mathbf{h}_s(\mathbf{x}_s, \mathbf{x}_{batch}, c_s,  t_f,  \mathbf{x}_{prod} ) = \mathbf{0} \\
            && \text{ } &&& \text{}  \mathbf{g}_s(\mathbf{x}_s, \mathbf{x}_{batch},  t_f,  \mathbf{x}_{prod} ) \leq \mathbf{0} \\
    \end{aligned}
\end{equation}

 In the \textit{approx tri} evaluations, $\mathcal{SC}(\cdot)$ takes the form of (\ref{eq: approx tri eval}) and approximates the optimal control response in (\ref{eq: bi eval}) for each scheduling problem by embedding the predicted optimal control responses as constraints into the scheduling via $\hat{c}_{ctrl}, \hat{t}^*_f = \mathcal{CTRL}( \mathbf{x}_{batch} )$. In each scheduling problem, we enumerate all 7 events by 2 machines by 5 product combinations of optimal control surrogates. The surrogates are trained as follows: we query the optimal costs and processing times corresponding to 10 uniformly distributed samples between 0 and the upper batch target limit for each equipment-production combination. We then train separate artificial neural networks with 2 hidden layers of 5 nodes each on the equipment-product datasets. No hyperparameter tuning is necessary at this point since this simple default architecture gives a good fit on the well-behaved optimal control response and remains tractable after embedding into the scheduling as discussed in Section \ref{sec: Surrogate accuracy}.

\begin{equation}
\label{eq: approx tri eval}
    \begin{aligned}
        & \mathcal{SC}( \mathbf{x}_{prod} ) =  &&  \underset{\substack{ \mathbf{x}_s,  \mathbf{x}_{batch},  c_s, \\   \hat{c}_{ctrl}, \hat{t}_f^*}}{ \arg \min.} 
            && f_s(\mathbf{x}_s, \mathbf{x}_{batch},  \hat{t}_f^*,  \mathbf{x}_{prod} ) \\
            & \text{ } && \text{s.t.} 
            && \mathbf{h}_s(\mathbf{x}_s, \mathbf{x}_{batch}, c_s,  \hat{t}_f^*,  \mathbf{x}_{prod} ) = \mathbf{0} \\
            && \text{ } &&& \text{}  \mathbf{g}_s(\mathbf{x}_s, \mathbf{x}_{batch},  \hat{t}_f^*,  \mathbf{x}_{prod} ) \leq \mathbf{0} \\
            && \text{ } &&&  \hat{c}_{ctrl}, \hat{t}^*_f = \mathcal{CTRL}( \mathbf{x}_{batch} )
    \end{aligned}
\end{equation}

In the \textit{tri} evaluations, two sequential optimization layers are solved within $\mathcal{SC}(\cdot)$: first, all 12 scheduling problems (\ref{eq: bi eval}) are solved in parallel. Then, the optimal batch responses $\mathbf{x}_{batch}^*$ for each scheduling problem as determined in (\ref{eq: bi eval}) feed into the 7 events by 2 machines optimal control problems as shown in (\ref{eq: tri eval}). The output of the control problem then returns the actual optimal control cost and processing times for feedback into the scheduling and planning level. In our case, the tri evaluations happen to be cheaper than the approx tri evaluations as discussed in Section \ref{sec: res meth2}. However, this was not known \textit{a priori} and is not expected to be the case when the size of the control problem increases. As a result, we use the tri evaluation only to check the solution quality in the first two methodologies and for minimal fine-tuning in the third.

\begin{equation}
    \label{eq: tri eval}
    \begin{aligned}[t]
        & \underset{\mathbf{x}_c, c_{ctrl}, t_f}{ \arg \min.} 
        && f_c(\mathbf{x}_c, t_f ,  \mathbf{x}^*_{batch} ) \\
        & \text{s.t.} 
        && \mathbf{h}_c(\mathbf{x}_c,  t_f, c_{ctrl},  \mathbf{x}^*_{batch} ) = \mathbf{0} \\
        &&& \mathbf{g}_c(\mathbf{x}_c, t_f,  \mathbf{x}^*_{batch} ) \leq \mathbf{0} \\
    \end{aligned} 
\end{equation}

Key to the tractable integration using DFO is parallelization in the scheduling layer. In our case, we see that the scheduling-level capacity for P1, P2, P3, and P4 production (Figure \ref{fig:iPSC}A) at any given time step is coupled to the planning level and other scheduling problems only via the inventory levels of I22. In the planning level, we are only interested in the changeover costs of the scheduling layer. As such, rather than solving the scheduling problems sequentially in time, we fix I22 inventory in the scheduling problems and solve all 12 scheduling problems in parallel leading to significant computational savings. In doing so, we rely on the planning-level safe storage constraints to ensure that we have enough I22 inventory for the production of any given planning target. Parallelization at the control level is in principle also possible, but in our case limited by the availability of compute nodes.

\subsection{Methodology 1: Derivative-free optimization} \label{sec: meth1}

The first methodology involves using DFO to find $\mathbf{x}_{prod}, \mathbf{x}_{transp}, \text{ and } \mathbf{x}_{sales}$ that optimize the solution metrics in (\ref{eq: metrics}). However, all solution metrics apart from planning-only (\ref{eq: metrics1}) are considered expensive as they call at least one additional level of scheduling optimization problems in each evaluation. As such, the DFO budget is severely limited to a couple of hundred evaluations given the high dimensionality of 178 DFO variables. This prevents even the most exploitative solvers from making meaningful progress. To mitigate this issue, we hot-start the DFO search by using the Pyomo \citep{pyomo} solution of the planning-only optimization problem (\ref{eq: planning-only opt}) as an initial guess. 

\begin{equation}
    \label{eq: planning-only opt}
    \begin{aligned}
            & \underset{\mathbf{x}_p, \mathbf{x}_{prod}}{ \text{min.}}
            && f_p(\mathbf{x}_p, \mathbf{x}_{prod})  \\
            & \text{s.t.}
            && \mathbf{h}_p(\mathbf{x}_p, \mathbf{x}_{prod}) = \mathbf{0} \\
            &&& \mathbf{g}_p(\mathbf{x}_p, \mathbf{x}_{prod}) \leq \mathbf{0} \\
    \end{aligned}
\end{equation}
Then, we sequentially optimize the bi (\ref{eq: metrics2}) and approx tri metrics (\ref{eq: metrics3}) using 500 evaluations each, where the solution of one feeds into the other one as an initial guess. 
We assume \textit{a priori} that the tri metric (\ref{eq: metrics4}) is intractable, and as such only use it to check solution quality. In principle, instead of this cascading approach, we could optimize the solution metric of interest from the start in the same solution time, which we discuss in Sections \ref{sec: res meth1} and \ref{sec: res meth3}.

Any DFO solver (e.g., Bayesian optimization, decision trees, direct methods) can be used for this task. We choose the state-of-the-art exploitative trust region optimizer Py-BOBYQA \citep{pybobyqa}, since the limited budget only allows for a fine-tuning of the solution around the initial guess rather than an extensive exploration of the solution space.

\subsection{Methodology 2: Optimality surrogates} \label{sec: meth2}

In Section \ref{sec: DFO init} we describe how we can use optimality surrogates to map the optimal control response as a function of the batch targets $\hat{c}_{ctrl}, \hat{t}^*_f = \mathcal{CTRL}( \mathbf{x}_{batch} )$. Similarly, for each planning step, we can construct optimality surrogates that return the optimal scheduling and control costs at a planning step as a function of the 4 production planning targets $\hat{c}^*_{s}, \hat{c}^*_c = [\mathcal{SCH-CTRL}]( \mathbf{x}_{prod} )$. The neural network $[\mathcal{SCH-CTRL}](\cdot)$ is trained to map the output of the (integrated) scheduling-control optimization problem$\mathcal{SC}(\cdot)$, which can take the form of the scheduling-control optimization problems (\ref{eq: bi eval})-(\ref{eq: tri eval}), associated with the bi, approx tri, and tri evaluation metrics respectively. As such, we could approximate the solution to the tri-level problem (\ref{eq: tri-level}) by solving the single-level problem (\ref{eq: surr upon surr}).

\begin{equation}
    \label{eq: surr upon surr}
    \begin{aligned}
            & \underset{\substack{\mathbf{x}_p, \mathbf{x}_{prod}, \\   \hat{c}_s^* \hat{c}^*_c}}{ \text{min.}}
            && f_p(\mathbf{x}_p, \mathbf{x}_{prod}) + \hat{c}^*_s +  \hat{c}^*_c  \color{blue}\\
            & \text{s.t.}
            && \mathbf{h}_p(\mathbf{x}_p, \mathbf{x}_{prod}) = \mathbf{0} \\
            &&& \mathbf{g}_p(\mathbf{x}_p, \mathbf{x}_{prod}) \leq \mathbf{0} \\
            &&&  \hat{c}_{s,t}^*, \hat{c}^*_{c,t}  =   [\mathcal{SCH-CTRL}]_t (\mathbf{x}_{prod}), \quad t=1, \dots, 12
    \end{aligned}
\end{equation}

Before training, we sample 1,000 combinations of the planning targets $\mathbf{x}_{prod}$ where the upper bound is obtained from the planning-level resource limit for each product. Each sample involves the solution of the integrated scheduling-control optimization problem $\mathcal{SC}(\cdot)$ as formulated in the scheduling-only (\ref{eq: bi eval}) or approximate scheduling-control (\ref{eq: approx tri eval}) formulations  associated with the bi and approx tri metrics respectively.  We use log-sampling within these bounds since uniform sampling oversamples infeasible operations as in practice the targets of the 4 products is rarely close to the upper bound. In the results section, we also repeat this workflow with surrogates trained on the accurate scheduling-control optimization problems (\ref{eq: tri eval}) to investigate the effect of inaccurate control information. We add fixed penalty terms to the planning objective if the scheduling (or control) problems become infeasible. Alternatively, adding soft penalties between the suggested and closest feasible set of planning targets to the DFO objective might lead to smoother surrogates.

Deciding on the architecture of the integrated scheduling-control surrogate is more difficult than in mapping the optimal control surrogate. The bigger the neural network, the more capable it is (in principle) of accurately capturing the optimal response function, but the more computationally expensive the planning formulation becomes after embedding of the surrogate. To navigate this accuracy-tractability trade-off, we use multi-objective Bayesian optimization for 12 evaluations to find the number of nodes in the first and second hidden layer $n_1, n_2$ that best trade-off accuracy and tractability. Each Bayesian optimization evaluation includes: the training of a network with hidden layers $n_1, n_2$; its embedding into the planning layer; optimization of the planning after embedding; and a solution quality check on all metrics of (\ref{eq: metrics}). The optimization time are returned with the approx tri (\ref{eq: metrics3}) or tri (\ref{eq: metrics4}) solution metrics and fed to the Bayesian optimization to suggest the next iteration of hyperparameters.

\subsection{Methodology 3: Combining both approaches}
The third methodology consists of combining the first two approaches to bring out the best of both worlds: We train a collection of the most promising network architectures identified in the off-line hyperparameter tuning from Methodology 2. Then, we optimize the planning problems with the embedded surrogates in parallel and evaluate their solution quality using the tri evaluation (\ref{eq: metrics4}) as the most accurate metric we can afford. Finally, we use DFO (in parallel) on the best solution(s) from the surrogates following Methodology 1 to fine-tune the solution towards the tri-level optimum using the tri evaluation (\ref{eq: metrics4}) in the remaining online evaluation budget.

\subsection{Software and practical considerations}

We use Py-BOBYQA \citep{pybobyqa} as our DFO solver of choice for the DFO approach (Methodology 1)  and the OpenBOX implementation \citep{openbox} for the multi-objective Bayesian optimization using Gaussian processes in the hyperparameter tuning as part of the surrogate approach (Methodology 2). The rest of the `conventional' optimization formulations are implemented in Pyomo \citep{pyomo}. While the control instances of Equation (\ref{eq: tri eval}) are solved using Ipopt \citep{ipopt}, the optimality surrogate formulation (\ref{eq: surr upon surr}) and all formulations for the integrated scheduling-control optimization $\mathcal{SC}(\cdot)$ from (\ref{eq: bi eval}) to (\ref{eq: approx tri eval}) are solved using Gurobi \citep{gurobi}. We use PyTorch \citep{pytorch} to train the neural networks and embed them as constraints into Pyomo blocks using OMLT \citep{omlt}. We run the main DFO and surrogate scripts on the Imperial College Research Computing HPC service. Most scripts run on a maximum walltime of 8 hours, and we increase the number of compute nodes and memory as needed, usually within the range of 64-128 and 10-50gb. Notable exceptions include the DFO approach on the approx evaluation using 18 hours of walltime, 256 nodes, and 25 gb memory, or the surrogate hyperparameter tuning on the tri evaluation using 8 hours of walltime, 128 nodes, and 100 gb.

In the next section, we present the accuracy-tractability trade-offs involved in incorporating the control surrogates into the scheduling, which constitutes the basis of moving from bi to approx tri evaluations. We then present the solution time and accuracy results of Methodology 1 and Methodology 2. In both approaches we are assuming that we have access to tri evaluations, but that they are expensive, meaning that we can only use them for checking the solution quality. In Methodology 3, we then present the solution time and accuracy results and discuss the advantages of combining both approaches together with relevant practical considerations. Throughout this, we place special emphasis in our discussion on distinguishing between the different solution accuracy metrics, between on- and off-line solution times, and the implications of model inaccuracies.

\section{Results and Discussion} \label{sec: Results}
Methodologies 1 and 2 heavily rely on the approx tri metric evaluation (embedding control surrogates into scheduling) as a proxy for the tri  evaluation. In the next section, we show why optimal control surrogates are well-suited for mapping the control problem, and we discuss the trade-offs involved in moving from bi to approx tri evaluations.

\subsection{Integrating scheduling and control} \label{sec: Surrogate accuracy}

\begin{figure}[htp!] 
    \centering
    \includegraphics[scale=0.5]{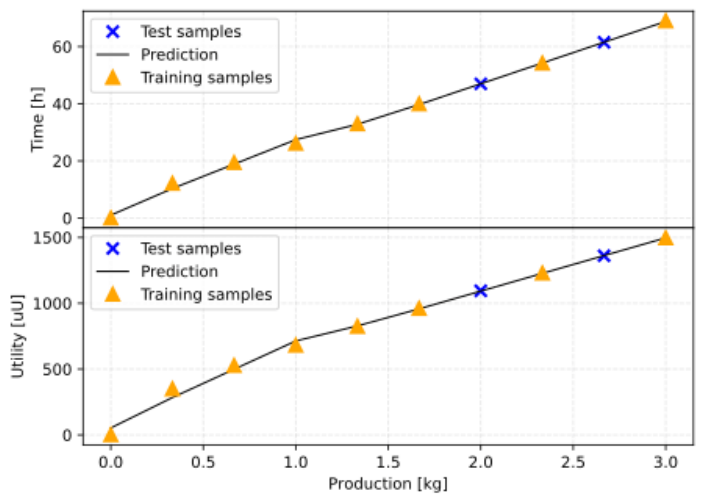}
    \caption{Predicted and actual optimal batch processing time and energy requirements corresponding to 10 uniformly sampled training and test batch targets.}
    \label{fig:Surrogate_Accuracy}
\end{figure}

Figure \ref{fig:Surrogate_Accuracy} gives an example of the predicted and actual optimal processing time and energy cost samples of the optimal control training set for one of the products. The mapping for both predictions is well-behaved, almost linear with a slight concave curvature. This mapping could be approximated reasonably well with linear or piecewise-affine linear regression techniques. 
This goes to show that even though the dynamics of the optimal control problem might be nonlinear, the relationship between the connecting variables in the integration of scheduling and control can be captured reasonably well with simple surrogates. In our case, a small neural network with 5 neurons for the two hidden layers each is sufficient. 

Table \ref{table:Ctrl_surr_time} summarizes the on- and off-line optimization times involved in the solution of the scheduling-only problem (\ref{eq: bi eval}) and the scheduling problem with embedded optimal control surrogates (\ref{eq: approx tri eval}). Solving these problems constitutes a key part in obtaining the bi and approx tri evaluations (\ref{eq: metrics2}) and (\ref{eq: metrics3}) respectively. The online optimization time increases from 10 seconds to over 4 minutes. Even though the surrogates are small, the approx evaluation requires the enumeration of each one of the 5 products by 7 events by 2 machines optimal batch production surrogates. Each of the 70 surrogates introduces discrete variables that lead to an explosion in solution time. While the online optimization time increases significantly after integration of the surrogates, the off-line control problem sampling and training of the neural network remain negligible. Training the neural network for 1,000 epochs takes at most a minute. Given the well-behaved nature of the mapping, we only require 10 control samples for each machine-product assignment, adding another one-off 5 minutes that are negligible in the overall workflow.

\begin{table}[htp!] 
 \centering
 \caption{Online optimization, control surrogate sampling and surrogate training times involved in the solution of the scheduling-only and scheduling with embedded control surrogates problems}
\begin{tabular}{c|ccc|}
\cline{2-4}
\multicolumn{1}{l|}{}                          & \multicolumn{3}{c|}{Time}                                                    \\ \hline
\multicolumn{1}{|c|}{}                         & \multicolumn{1}{c|}{Optimization} & \multicolumn{1}{c|}{Control sampling (N=10x5)} & Training \\ \hline
\multicolumn{1}{|l|}{Scheduling only}          & \multicolumn{1}{c|}{10 sec}        & \multicolumn{1}{c|}{/}    & /    \\ \hline
\multicolumn{1}{|c|}{Scheduling-control surrogate}       & \multicolumn{1}{c|}{250 sec}        & \multicolumn{1}{c|}{5 min}    & 1 min    \\ \hline
\end{tabular}
 \label{table:Ctrl_surr_time}
 \end{table}


\subsection{Methodology 1: Integration into planning via DFO} \label{sec: res meth1}

Optimizing the scheduling-only problem or scheduling problem with embedded surrogates as presented in the previous section constitutes the bottleneck in the bi and approx tri evaluations. Figure \ref{fig:DFO_approach} illustrates the solution accuracy obtained and time required to perform DFO on the bi and then on the approx tri metrics for 500 evaluations each. 

\subsubsection{Solution accuracy}

Figure \ref{fig:DFO_eval} shows the different solution metrics obtained in the DFO. We start with the bi DFO (DFO on the bi metric) using the solution obtained from the optimization of the planning-only problem in Pyomo. This planning-only solution is infeasible for all evaluation metrics accounting for lower level information.

The planning-only solution then feeds into the bi DFO where the bi metric is improved from around 6 to -3 in 500 evaluations. Then, the bi DFO solution feeds into another round of 500 approx tri DFO evaluations, improving the approx metric from the bi DFO solution from about -2 to -3. This goes to show that DFO reliably improves the solution according to the evaluation metric used in the objective. Yet, there is no guarantee that an increase in cheaper bi/approx evaluations comes with an improvement in more accurate evaluations that are not used as the objective function of the DFO. This explains why the tri evaluation shows no improvement - all 3 tri evaluations remain infeasible in some scheduling and/or control problems. 

Realistically, given the expense of the scheduling-level optimization problems and the high dimensionality in the planning-level variables, we only have enough budget to fine-tune the solution obtained in the initial guess. This is confirmed in \ref{sec: Solution profiles} showing the planning profile associated with the DFO solutions. However, this fine-tuning can render the planning-only solution feasible with respect to more accurate metrics, resulting in large increases in solution quality. While we use fixed costs to penalize infeasibilities in the scheduling and control level, DFO would make more consistent progress if soft penalty violations would be used between the suggested and nearest feasible set of planning and batch targets. 

\subsubsection{Solution time}

\begin{figure}[htp!] 
    \centering
    \begin{subfigure}[b]{0.475\textwidth}
         \centering
             \includegraphics[width=\textwidth]{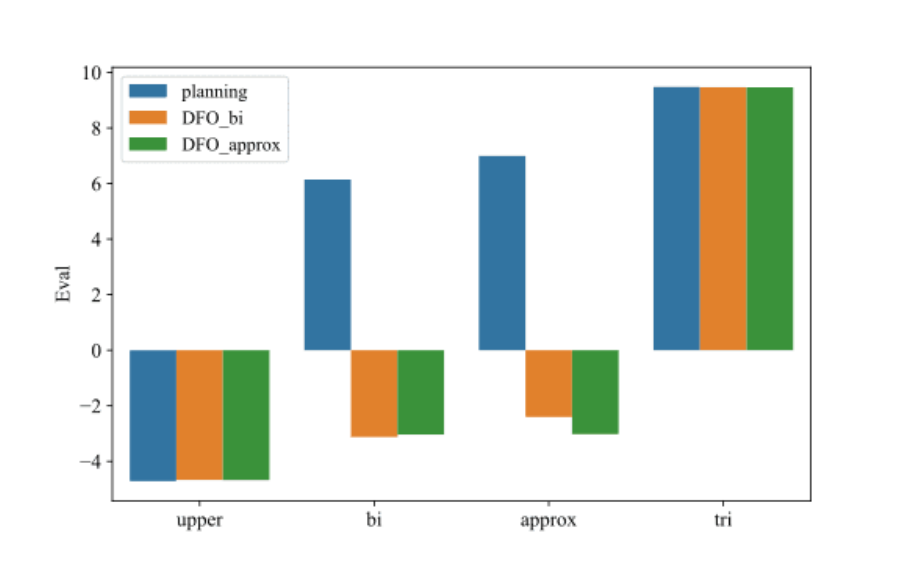}
         \caption{Evaluation achieved on the upper, bi, approx tri, and tri solution metrics}
         \label{fig:DFO_eval}
     \end{subfigure}
     \begin{subfigure}[b]{0.475\textwidth}
         \centering
             \includegraphics[width=\textwidth]{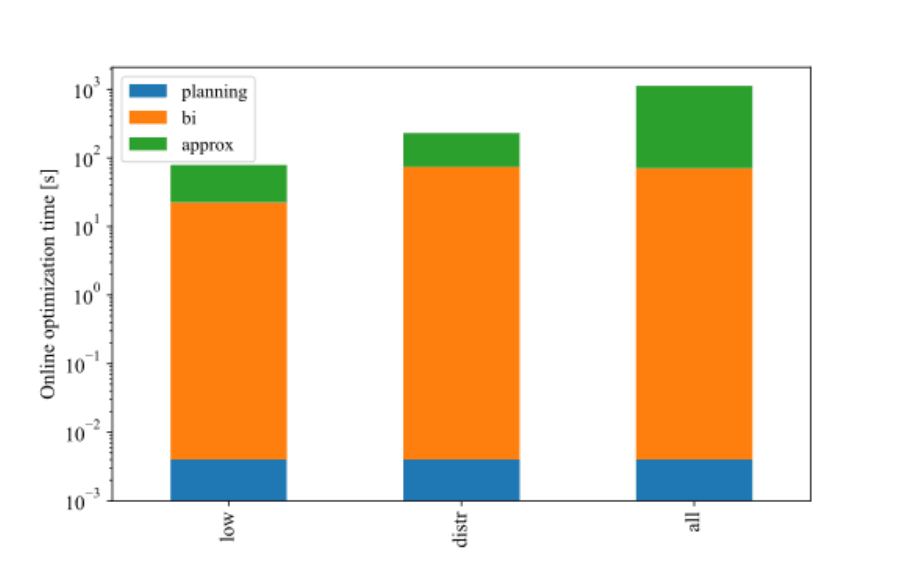}
         \caption{Cumulative online optimization time on the low, distr, and full configurations}
         \label{fig:DFO_time}
     \end{subfigure}
     \caption{Solution accuracy and time for the planning-only instance, the DFO approach on the bi metric, and the DFO approach on the approx tri metric}
     \label{fig:DFO_approach}
\end{figure}

Figure \ref{fig:DFO_time} shows the solution times involved in finding the planning-only initial guess and in subsequently performing 500 DFO evaluations on the bi and approx metric each on three configurations of the case study: in \textit{low}, we only solve a single average scheduling problem for all planning steps wherein we account for the optimal control of only a single product; in \textit{distr}, we solve in parallel a separate scheduling problem for each planning step accounting for the optimal control of only a single product; \textit{all} refers to the full case study as described in Section \ref{sec: Case study}, where at each planning step, we solve in parallel all scheduling problem where the control of all 5 products is performed optimally. Comparing the `low' and `distr' configurations allows to investigate the effect of parallelization while comparing the `distr' and `all' configurations allows to investigate the effect of increasing the number of control problems per scheduling instance.

We see that the planning-only Pyomo instance returns a solution in milliseconds. When going from the `low' to the `distr' configuration, the bi and approx DFO times only increase by a factor of 3-4 from around 20 to 80 and from around 60 to 200 seconds. This highlights the benefits of parallelization: rather than seeing an increase in computational time proportional to the number of planning steps (12), we only see it increase by a factor of around 4. In the best case scenario, this parallelization would keep the solution time constant. In practice however, each evaluation call is limited by the time it takes for the slowest scheduling instance to be solved. 
Increasing the number of optimal control problems when moving from `distr' to `all' only increases the approx DFO time. Increasing the number of control problems by 5 increases the solution time by a factor of 10, from 100 to almost 1,000 seconds. This is to be expected given the discrete nature of the variables introduced in the surrogates. In Table \ref{table:Ctrl_surr_time}, we claim that the optimization instances in a single bi or approx evaluation take 10 and 250 seconds respectively. However, Figure \ref{fig:DFO_time} shows that 500 evaluations take around 100 and 1,000 seconds respectively. While this may seem contradictory, this means that most of the DFO samples are quickly determined to be infeasible. 

\subsubsection{Cascading approach}

The question arises if rather than use a cascading approach of feeding the DFO solutions of less accurate metrics as initial guesses into the DFO of more accurate metrics, we can obtain better solutions by using the available budget to have more samples to optimize the expensive metric of interest directly. The answer to this is case-dependent and is influenced by the extent to which the metrics corresponding to the different levels of integration are correlated. In our case, the cascading approach is justified since the DFO on the bi metric is cheaper by an order of magnitude and still makes significant progress on the approx evaluation. However, if the ultimate goal is to optimize the tri evaluation, then neither the optimization of the bi nor the approx metric make significant progress. We continue this discussion in Section \ref{sec: res meth3}.

\subsection{Methodology 2: Integration into planning via surrogates} \label{sec: res meth2}

In Methodology 2, rather than perform DFO on the bi or approx evaluations, we train optimality surrogates as a function of the planning targets to map the scheduling and control costs from the scheduling-only  (\ref{eq: bi eval}) and approximate scheduling-control problems (\ref{eq: approx tri eval}). Since these instances constitute a major part of the bi and approx tri evaluations, we refer to the surrogates trained on them as bi or approx surrogates respectively. In Figure \ref{fig:Surrogate_approach}, for the bi and approx surrogates, we present three architectures each that return a Pareto solution between solution time and best approx evaluation after optimization of the planning level with the embedded surrogate.

\subsubsection{Solution accuracy}

\begin{figure}[htp!] 
    \centering
    \begin{subfigure}[b]{0.475\textwidth}
         \centering
             \includegraphics[width=\textwidth]{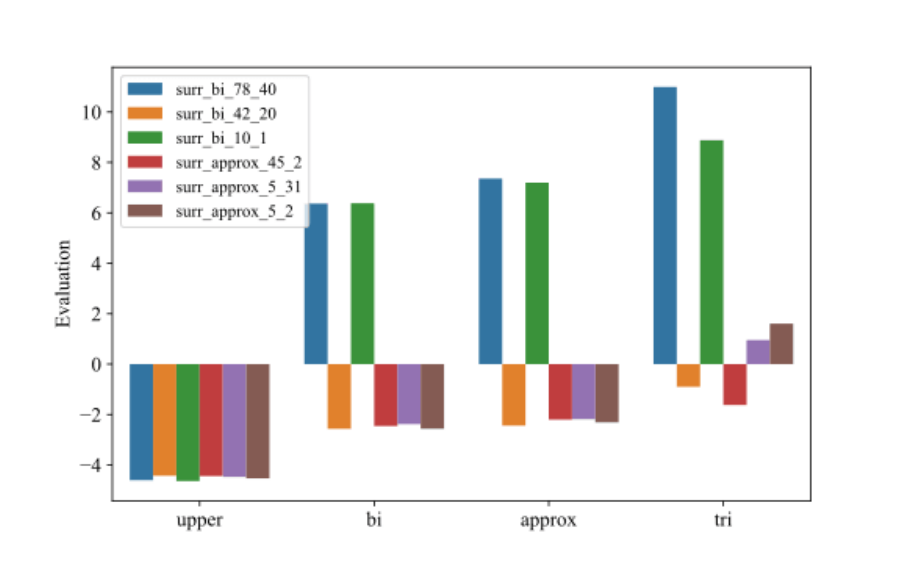}
         \caption{Evaluation achieved on the upper, bi, approx, and tri solution metrics}
         \label{fig:Surr_eval}
     \end{subfigure}
     \begin{subfigure}[b]{0.475\textwidth}
         \centering
             \includegraphics[width=\textwidth]{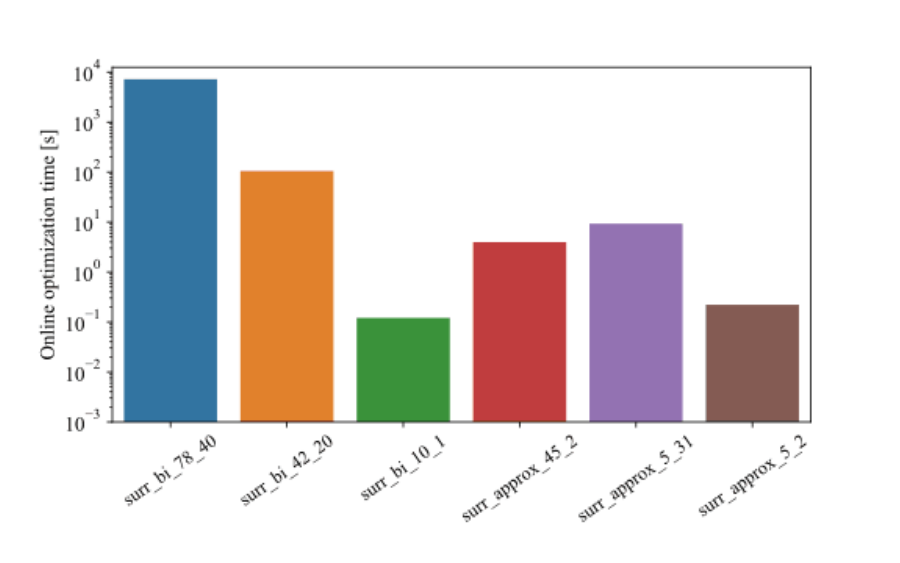}
         \caption{Online optimization time for the planning problem after embedding of the surrogates}
         \label{fig:Surr_time}
     \end{subfigure}
     \caption{Three solution accuracy and time Pareto solutions each of surrogates trained on the bi and approx tri metrics with the two numbers in the label referring to the number of nodes in the two hidden layers}
     \label{fig:Surrogate_approach}
\end{figure}

Figure \ref{fig:Surr_eval} shows that the solutions obtained of the approx surrogates are more consistent and achieve better solution quality than those obtained by the bi surrogates. This suggests that using more informative evaluations to build the surrogate translates into better solution accuracy. As opposed to the DFO where using the bi evaluation metric as the objective significantly improves the bi evaluation, the bi surrogates fail to consistently find a good bi evaluation. This might be because of the way the surrogates are trained. In the hyperparameter tuning, both the bi and approx surrogate architectures are optimized based on the approx metric evaluation, as it was expected \textit{a priori} to be the most accurate computationally affordable metric available. The bi surrogates, only trained on scheduling-only problems used in the bi evaluations, might miss valuable information that relates surrogate architecture, training accuracy, and optimization solution quality, leading to poorer, less consistent performance.

However, it is surprising to see that the solution check on the tri metric evaluations corresponding to the approx surrogates are significantly better than the ones obtained in the DFO. This suggests a different explanation for why the surrogate approach performs surprisingly well. Rather than the surrogate approach performing well \textit{despite} the inevitable small model inaccuracies, the surrogate approach might perform well \textit{because} of them. In the DFO solutions, we are starting off from the planning-only solution that is optimal in the planning-level formulation only. This planning-only solution is probably active in some constraints that quickly become infeasible when more accurate information from lower-level scheduling or control becomes available. Small inaccuracies in the surrogates might give a suboptimal solution that is further away from the constraints and as such more `robust' to uncertainties in the lower level. This hypothesis is supported by the planning solution profiles of the surrogate approaches being significantly different to the profiles observed for the centralized and DFO solutions as seen in \ref{sec: Solution profiles}. It is in principle possible that model inaccuracies might be biased towards giving solutions in the infeasible space, but these models should be filtered out by the hyperparameter tuning process being biased towards models giving `safer' solutions.

\subsubsection{Consistency in solution accuracy}

\begin{figure}[htp!] 
    \centering
     \includegraphics[width=0.475\textwidth]{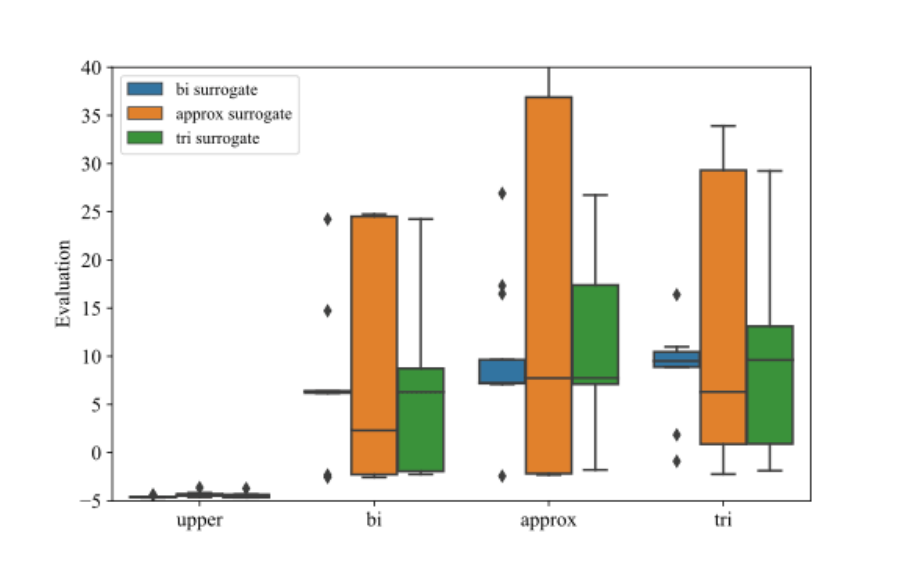}
     \caption{Box-and-whisker plot of the upper, bi, approx tri, and tri solution metric evaluations associated with all 12 samples of the hyperparameter tuning for the bi, approx tri, and tri surrogates. The box represents the quartiles while the whiskers demonstrate variability outside the box plots. Outliers are shown as isolated data points outside of the whiskers.}
     \label{fig:Surrogate_var}
\end{figure}

The major risk of the surrogate approach lies in selecting an architecture or training a specific model that leads to poor performance. To determine the effect of the model architecture on performance, Figure \ref{fig:Surrogate_var} displays the solution accuracy associated with the 12 model architectures queried in the hyperparameter tuning. Additionally, we investigate how much of the variability in bi and approx tri surrogate model performance is due to them not being trained on accurate control information. To this end, we follow the same workflow to train tri surrogates on accurate scheduling-control instances as used in the tri metric (\ref{eq: metrics4}). 

Overall, the best of the 12 samples display similar performance among the three surrogate types on all solution metrics. In general, we see that the bi and tri surrogates display a similar median and variation in solution accuracy, with the bi surrogates being slightly more consistent in the tri metric evaluation. While the approx surrogates display a better median evaluation in the bi and tri metric evaluation, they also display more variability in solution performance. This is surprising to see given Figure \ref{fig:Surr_eval} where the best approx tri surrogate Pareto solutions seem to be more consistent than their bi surrogate counterparts. This can be explained however by the approx surrogates being not just subject to model inaccuracy in the scheduling-control surrogates but themselves already trained on imperfect optimal control surrogates. 

The question still remains as to the tri surrogates' inability to clearly outperform the other surrogates. One possible explanation could be that accounting for the optimal control makes the solution space too ill-behaved to be accurately captioned by the networks.
Overall, the surrogates display too similar performance to argue that one type outperforms the other. Ideally, in the hyperparameter tuning, we would also account for the sensitivity of the surrogate performance with respect to different training configurations or slightly updated/modified data, which calls for additional work on model validation or robust/adversarial selection approaches. 

\subsubsection{Solution time}

Figure \ref{fig:Surr_time} shows the online optimization time required to solve the embedded problem as shown in Formulation (\ref{eq: surr upon surr}). The online optimization times of the bi and approx tri surrogates range widely based on the size of the networks from 100 milliseconds to 2 hours. It also seems like the type of surrogate matters less than the number of hidden nodes. For example, the bi surrogate with hidden layers of sizes 10 and 1 leads to a solution time comparable to that of the approx surrogate with sizes 5 and 2. The biggest surrogate (78 and 40 nodes) times out after 2 hours of solution time and displays the worst of the Pareto solutions in terms of accuracy. For the bi surrogates there seems to be a trade-off between using bigger and more accurate surrogates leading to better solutions at the expense of high computational costs as the best solution is attributed to the medium-sized surrogate of 42 and 20 nodes. However, for the approx surrogates, there does not seem to be much difference in solution quality found, since the best solution is found using the smallest network. This makes intuitive sense: by constructing an approximation of an approximation, approx surrogates might present similar, smoothed response surfaces across architectures.

While Figure \ref{fig:Surr_time} suggests that the online optimization time depends primarily on the surrogate size rather than the type of evaluation it is trained on, Table \ref{table: surr} shows that the difference in solution time between the different types of surrogates is shifted off-line to the sampling, training, and hyperparameter tuning steps. The training times are essentially negligible with each surrogate being trained in a maximum of 90 seconds. The hyperparameter tuning time can be adjusted based on the available budget. In our case, for each surrogate type, we limit the online optimization time per embedded problem to 2 hours and the hyperparameter tuning time involving 12 surrogates to 8 hours. Essentially, the bottleneck in moving from bi to approx surrogates is in the sampling step. We sample 1,000 bi, approx, and tri evaluation samples, which takes between 8 and 24 hours. As discussed in the DFO section, the 8 and 24 hours for the bi and approx evaluation samples are considerably less than what would be expected by extrapolating the evaluation times shown in Table \ref{table:Ctrl_surr_time} by the number of samples. This is expected since only around half of the samples turn out to be feasible. 

It is surprising to see that the tri sampling time is shorter than the approx sampling time.
It turns out that solving the planning problems, the scheduling problems, and then the 168 (2 machines by 7 events by 12 planning steps) optimal control problems in sequence is 5 times faster on average than solving the scheduling problems with 70 optimal control surrogates (7 events by 2 machines by 5 productions) embedded.
In other words, it is not trivial to see that embedding all possible optimal control combinations into the scheduling displays worse scaling than evaluating the tri metric by solving the levels hierarchically.
We expect that as the optimal control problems become more expensive, the tri evaluation becomes more expensive than the approx evaluation. In hindsight, for the case study at hand, we conclude that since the tri evaluation is faster than the approx tri evaluation, the performance of using approx tri evaluations in the DFO or their corresponding approximate scheduling-control problem (\ref{eq: approx tri eval}) for surrogate training is not worth the compute.

\begin{table}[htp!] 
\centering
\caption{On-line optimization, and off-line hyperparameter tuning, training, and sampling times associated with the construction of surrogates on the scheduling-only, approximate scheduling-control, and accurate scheduling-control}
\label{table: surr}
\begin{tabular}{cc|cccc|}
\cline{3-6}
                                                                                               &                                                              & \multicolumn{4}{c|}{Time}                                                                                                                                                                                                     \\ \cline{2-6} 
\multicolumn{1}{c|}{}                                                                          & \begin{tabular}[c]{@{}c@{}}Related\\ evaluation\end{tabular} & \multicolumn{1}{c|}{Online}      & \multicolumn{1}{c|}{\begin{tabular}[c]{@{}c@{}}Hyperparameter\\ tuning\end{tabular}} & \multicolumn{1}{c|}{Training}        & \begin{tabular}[c]{@{}c@{}}Sampling\\ (N=1,000)\end{tabular} \\ \hline
\multicolumn{1}{|c|}{Scheduling-only}                                                          & Bi                                                           & \multicolumn{1}{c|}{0.1min-2hrs} & \multicolumn{1}{c|}{8 hrs}                                                           & \multicolumn{1}{c|}{up to 12x90 sec} & 8 hours                                                      \\ \hline
\multicolumn{1}{|c|}{\begin{tabular}[c]{@{}c@{}}Approximate\\ scheduling-control\end{tabular}} & \begin{tabular}[c]{@{}c@{}}Approx \\ tri\end{tabular}        & \multicolumn{1}{c|}{0.1min-2hrs} & \multicolumn{1}{c|}{8 hrs}                                                           & \multicolumn{1}{c|}{up to 12x90 sec} & 24 hours                                                     \\ \hline
\multicolumn{1}{|c|}{Scheduling-control}                                                       & Tri                                                          & \multicolumn{1}{c|}{0.1min-2hrs} & \multicolumn{1}{c|}{8 hrs}                                                           & \multicolumn{1}{c|}{up to 12x90 sec} & 8 hours                                                      \\ \hline
\end{tabular}
\end{table}

\subsection{Comparing the DFO and surrogate approaches}
In the above two sections, we highlight the accuracy-tractability trade-offs associated with both approaches.
DFO does not introduce additional inaccuracies if the lower levels are solved to global optimality.
As such, DFO guarantees a solution at least as good as the initial guess on the evaluation metric used as its objective.
The main drawback of DFO is its poor scaling in the number of planning-level variables and in the computational expense of the lower-level problems.
Realistically, given these budget restrictions, we can only expect to fine-tune the solutions around the initial guess.

Generally, the surrogate approach scales better at the potential expense of solution quality. A main advantage of the surrogate approach is that tractability and accuracy can be easily adapted to the problem at hand by using Bayesian optimization (which ironically counts as DFO too).
Empirically, we see that even the optima obtained with smaller, more tractable, and potentially less accurate models can display tri-level metric solutions that outperform any of the solutions obtained via DFO. However, by using surrogates, we forfeit any guarantee in the solution quality obtained. While surrogate techniques promise on-line computational savings, we have to invest our efforts off-line into the model selection and validation process.
In the next section, we synthesize these findings into a proposed method that leverages the advantages of both methods.

\subsection{Methodology 3: Bringing out the best in both approaches} \label{sec: res meth3}

Figure \ref{fig:Best_approach} shows the improvement in solution quality and increase in computational expense incurred by running DFO using the tri evaluation metric on the planning-only, DFO and surrogate solutions.

\subsubsection{Solution accuracy}

Figure \ref{fig:Best_eval} shows that running 500-1,000 DFO evaluations using the tri metric on the previous DFO and surrogate solutions significantly improves their tri evaluation.
The best tri evaluations are found by performing DFO on the bi and tri surrogate solutions for a standard 500 evaluations, improving the tri metric from about -1 and -2 to around -2.5 respectively. We use the same standard budget on the approx DFO solution, improving the solution from about 9.5 to 8.5. However, by running DFO on the planning-only and bi DFO solutions using 1,000 and 800 evaluations, we find a solution that is just as good in less time than is required to perform the approx DFO in the first place. Performing tri DFO on the bi DFO solution even improves the solution to 7.5. This links back to the previous discussion in Section \ref{sec: res meth1} of whether it is favourable to perform DFO on increasing levels of metric complexity in a cascading manner or if we should use the available budget to only optimize the targeted metric from the start. In our case, we find that we should use the cascading approach but avoid the more expensive approx evaluation in favour of the tri evaluation.

\subsubsection{Solution time}

\begin{figure}[htp!] 
    \centering
    \begin{subfigure}[b]{0.475\textwidth}
         \centering
             \includegraphics[width=\textwidth]{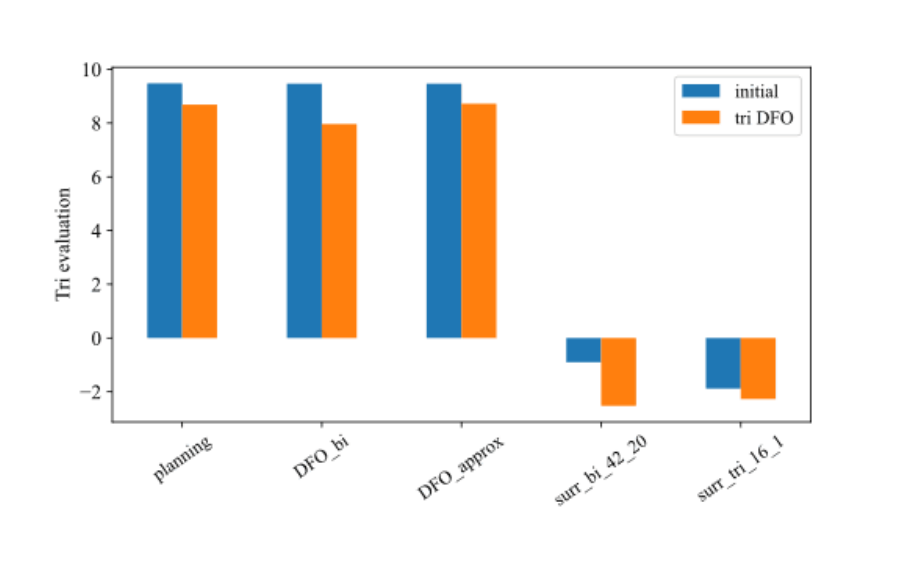}
         \caption{Tri evaluation of the initial guess and subsequent tri DFO solution}
         \label{fig:Best_eval}
     \end{subfigure}
     \begin{subfigure}[b]{0.475\textwidth}
         \centering
             \includegraphics[width=\textwidth]{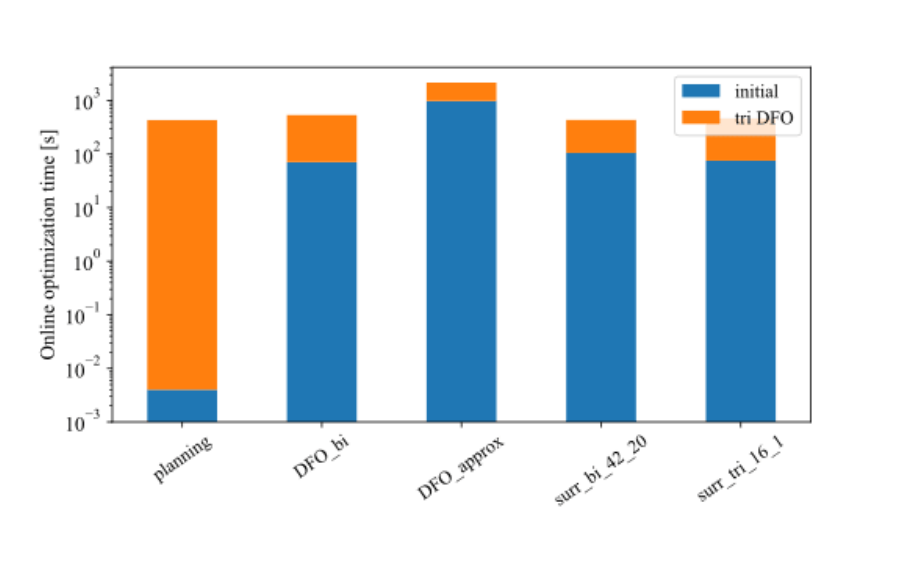}
         \caption{Cumulative online optimization time to find the initial guess and perform the tri DFO}
         \label{fig:Best_time}
     \end{subfigure}
     \caption{Increase in solution time and quality in using DFO to optimize the planning-only, DFO, and surrogate solutions on the tri evaluation metric}
     \label{fig:Best_approach}
\end{figure}

Figure \ref{fig:Best_time} compares the computational time required to fine-tune the solutions using tri DFO with the time it takes to obtain the initial guesses in the first place. 
500 tri DFO evaluations on the approx DFO and surrogate solutions take about 240 minutes, while the initial 500 approx DFO evaluations take about 1,000 minutes in the first place. 1,000 tri DFO evaluations take about 430 minutes starting from the planning-only solution to get to a final tri evaluation of 8.5. 800 tri DFO evaluations on the bi DFO solution take about 390 minutes to achieve a solution of 7.5. In both cases, the tri DFO constitutes the computational bottleneck compared to the 0.2 and 70 minutes needed to obtain the planning-only and bi DFO solutions in the first place. Either way, both approaches take significantly less time than the complete cascading approach including the 500 bi, approx and tri DFO evaluations taking over 20 hours cumulatively.
We note again that the average DFO evaluation time is skewed by quick infeasible runs.

While the bi DFO approach outperforms the other DFO approaches for a given time budget, we conclude that better solutions can be obtained in less time when hot-starting the tri DFO with the surrogate approach. In our case, the 500 tri DFO evaluations take about 2-3 times longer than the surrogate optimization. The relative solution time between the tri evaluations and surrogate optimization changes significantly with the number of infeasible samples in the DFO and the size of the surrogates. However, this ratio can be manipulated by the DFO budget and the maximum surrogate runtime.

\subsubsection{Practical considerations}

Optimality surrogates elegantly trade-off tractability and accuracy. However, their variability in solution performance is concerning. 
In order to mitigate these drawbacks, we suggest finding safety in numbers and choosing multiple surrogates - either the same model architectures trained on slightly different datasets, or different architectures potentially even trained on different kinds of evaluations. We can then embed these surrogates separately into different planning instances and optimize the resulting formulation in parallel in the online optimization step. All solutions are then checked on the most accurate evaluation metric available. Finally, we perform DFO on the most accurate evaluation metric in the hopes of finding a better solution in the remaining time budget. At this stage, the solution quality is guaranteed to be as good as the best surrogate solution for a chosen evaluation metric. Alternatively, we can still fall back on the much cheaper planning-only solution, which represents the traditional `sequential' solution approach.

This leverages the available on- and off-line optimization budget to its fullest. We can leverage parallelization to mitigate model inaccuracy risks on-line. We also still retain the advantage of the surrogate approach of having shifted a lot of compute off-line. The amount of effort put into the offline sampling, training, hyperparameter tuning, and model validation step can be adjusted to the tolerance for inaccuracies in the problem at hand. The remainder of the available online solution time can then be used to its fullest to fine-tune the surrogate solutions via DFO. 

\section{Conclusion} \label{sec: Conclusion}
In this work we present how derivative-free optimization (DFO) and optimality surrogates can be used to optimize general bi-level formulations. We show how these techniques can be used to increase the solution quality of a multi-level, multi-site hierarchical planning-scheduling-control problem with respect to the conventional sequential solution approach. In tackling a task as ambitious as trying to solve these large-scale tri-level formulations, we need to use tools that improve the solution quality in as little time as possible. At the same time, we need to mitigate the risk of the solution becoming unusable through model inaccuracies. We synthesize our findings into a methodology that leverages the advantages of both techniques: We use derivative-free optimization (Bayesian optimization) in the off-line model selection process to efficiently trade-off solution accuracy and tractability. We then mitigate the risk of model inaccuracies in the surrogates by employing the surrogate approach in parallel on a combination of multiple model instances. After verifying the solution quality on the most accurate affordable solution metric, we can finally use DFO to fine-tune the best solution(s) in the remaining online optimization budget. 

While we can make general observations about the trade-offs involved in both methods, it should be noted that these are limited to the case study at hand. While the DFO approach is agnostic to the integrated scheduling-control formulation, the efficacy of the approach depends on the relative size of the three layers. Similarly, the surrogate approach is agnostic to the lower-level problems as long as the bottom-up complicating variables can be mapped as a function of the top-down complicating variables. Future work would investigate and stress-test different planning network configurations, different types of scheduling formulations, and different surrogate formulations in the control.
The surrogate technique in particular opens up many avenues of research. For instance, instead of using multi-objective Bayesian optimization, we could optimize only solution quality under time limit constraints. We should also investigate how to account for model consistency in the hyperparameter tuning, i.e. the variability of model architecture performance when trained on different data. 
Ultimately, our work demonstrates the level of maturity of data-driven techniques and their ease of integration into more conventional optimization frameworks. This prompts us to re-evaluate their potential in addressing more accurate hierarchical formulations, even in larger-scale, industrial case studies.

\paragraph{Acknowledgements}
D.v.d.B. gratefully acknowledges financial support from the EPSRC DTA studentship award and the Bansal bursary. We also acknowledge computational resources and support provided by the Imperial College Research Computing Service (http://doi.org/10.14469/hpc/2232)

\paragraph{Data Availability Statement}
All data necessary for replication of results can be found with the code under \href{https://github.com/OptiMaL-PSE-Lab/DD-Hierarchical}{https://github.com/OptiMaL-PSE-Lab/DD-Hierarchical}.

\listoffigures

\clearpage


\bibliography{biblio}
\clearpage

\appendix

\section{Optimization formulations} \label{sec: Appendix opt}

\subsection{Planning formulation} \label{sec: Appendix plan opt}

Formulation (\ref{eq: app Planning}) minimizes the planning objective (\ref{eq: app plan obj}) subject to constraints (\ref{eq: app plan mass}) to (\ref{eq: app plan safe}) . The first term includes the material production and storage costs over all sites $l \in L=\{S1,S2,S3\}$, materials $s \in S=\{I1, I2, I21, I22, P1,P2,P3,P4,P5,P6\}$, and timesteps $t \in T=\{1,\dots, 12\}$; The second term consists of the transport costs of I2 to sites S2 as I21 and to S3 as I23; The third term presents the income from selling the products $p \in I_p=\{P1, \dots, P6\}$. Variable $P_{i,l,t}$ is the production of material i in location l at time t produced at a cost $CP_i$;  Variable $S_{i,l,t}$ denotes the inventory of material i in location l at time t stored at a cost $CS_i$ with initial inventory $S0_{i,l}$; Variable $TP_{l,tc}$ is the transport of material I2 to sites S2 or S3 at the transport control step $tc \in Tc = \{1,\dots,5\}$ at a cost $CT_l$; Variable $SA_{p,t}$ is the sales of product p at time t at the price of $SP_p$.

\begin{subequations}
\label{eq: app Planning}
\begin{align}
     & \min.  
     && \sum_{i \in I, t \in T, l \in L} ({CP}_{i} P_{i, l, t} + {CS}_{i} S_{i,l,t}) + \sum_{l \in \{S2, S3\}, tc \in Tc} {CT}_{l} {TP}_{l,tc} - \sum_{p \in I_p,t} {SP}_{p} {SA}_{p,t} \label{eq: app plan obj}\\
     & \text{s.t.} 
     &&  S_{i,l,t} = S_{i,l,t-1} + P_{i,l,t} - {Cons}_{i,l,t} - {TR}_{out, i,l,t} + {TR}_{in, i,l,t}, \quad i \in I, l \in L, t \in T \setminus \{0\} \label{eq: app plan mass} \\
     &&&  S_{i,l,0} = {S0}_{i,l} + P_{i,l,0} - {Cons}_{i,l,0} - {TR}_{out, i,l,0}, \quad i \in I, l \in L \label{eq: app plan mass2}\\
     &&& {Cons}_{i,l,t} = 1.1 \sum_{i' \in I_{out}(i)} Q_{i',l} P_{i',l,t},  \quad i \in I, l \in L, t \in T \label{eq: app plan cons} \\
     &&& {TR}_{out, i,l,t} =  \mathbbm{1}_{i=I2,l=S1} \sum_{l' \in \{S2,S3\}} {TP}_{l,to\_tc(t)}, \quad i \in I, l \in L, t \in T  \label{eq: app plan trout}\\
     &&& {TR}_{in, i,l,t} =   \mathbbm{1}_{(i,l) \in \{(I21, S1),(I22,S2)\},t-{LT}_{l}\geq 0} \quad {TP}_{l, to\_tc(t-{TT}_{i,l})}, \quad i \in I, l \in L, t \in T \label{eq: app plan trin} \\
     &&& P_{i,l,t} \leq M_{big} X_{l,i}, \quad i \in I, l \in L, t \in T \label{eq: app plan prod2loc} \\
     &&&  U_{l,r} \sum_{p \in I_L(l)} Q_{p,l} P_{p,l,t} \leq A_{l,r}, \quad r \in R, l \in L, t \in T \label{eq: app plan res} \\
     &&& {SA}_{p,t} \leq D_{p,t}, \quad p \in I_p, t \in T \label{eq: app plan demand} \\
     &&&  S_{i,l,t} \geq S^*_p \text{ if } t > 4, \text{ else } S_{i,t} \geq 0.25*S^*_p , \quad i \in I, l \in L, t \in T \label{eq: app plan safe}
\end{align}
\end{subequations}

Constraints (\ref{eq: app plan mass}) and (\ref{eq: app plan mass2}) define the inventory mass balances, where variable ${Cons}_{i,l,t}$ is the consumption of material i at site l and timestep t and is related in (\ref{eq: app plan cons}) to the production of subsequent materials $I_{out}(i)$ via the conversion factor $Q_{i',l}$. Variables ${TR}_{out, i,l,t},{TR}_{in, i,l,t}$ define the removal and addition of material i at site l and timestep t via transport as defined in (\ref{eq: app plan trout}) and (\ref{eq: app plan trin}): ${TR}_{out, i,l,t}$ is only relevant to I2 at site 1 and sums up the transport of material I2 to sites 2 and 3 at time t, whereas ${TR}_{in, i,l,t}$ is only relevant to I21 and I22 at sites 2 and 3 if material I2 was transported from site 1 $LT_l$ (leadtime) steps away. $to\_Tc(\cdot)$ maps the 12 planning-level timesteps $t$ to the 5 transport timesteps $tc$. Parameter $X_{l,i}$ denotes if a material i is produced at site l and sets the production of non-relevant materials for a given site to zero in (\ref{eq: app plan prod2loc}) with the help of a big-M constant $M_{big}$. $A_{l,r}$ in (\ref{eq: app plan res}) denotes the upper limit on the resource utilization $U_{l,r}$ associated with all production at a certain site l for a given resource r at a given time. (\ref{eq: app plan demand}) upper bounds the sales of all products by their demand $D_{p,t}$. \ref{eq: app plan safe} constrains the inventories of materials $S_{i,l,t}$ to be above their safe storage $S^*_p$ at all timesteps, with slight relaxations in the first timesteps.

\subsection{Scheduling formulation} \label{sec: Appendix sch opt}

The optimal planning targets $P_{p,S2,t}$ for each product $p \in \{P1,P2,P3,P4\}$ of site S2 in the 12 planning timesteps are fed to a separate scheduling formulation as $P_s$. All scheduling formulations optimize the makespan $MS$ while penalizing the number of changeovers $Y_{i,i',n}$ required to achieve the planning targets as determined in the planning problem. We adapt \citet{scheduling} and rely on the state-task network approach with a common continuous-time representation for all units. The formulation also accounts for variable batch sizes, variable processing times and sequence-dependent changeover times. The formulations for the assignment constraints (\ref{eq: app Scheduling assignment}), duration, finish time, and time-matching constraints (\ref{eq: app Scheduling duration}), batch size constraints (\ref{eq: app Batch duration + MB}), and tightening constraints (\ref{eq: app Scheduling tightening}) are presented in great detail in \citet{scheduling}.

Objective:
\begin{subequations}
\label{eq: app Scheduling obj}
\begin{align}
     &  MS + \rho \sum_{i \in I, i' \in I, n \in N} Y_{i, i', n} 
\end{align}
\end{subequations}
Variable $MS$ denotes the makespan, and parameter $\rho$ penalizes the sum of changeovers $Y_{i,i',n}$ from task $i \in I=\{I23M1, P1M1,P2M1,P3M1,P4M1, I23M2, P1M2,P2M2,P3M2,P4M2\}$ to task $I' \in I$ at scheduling event $n=\{1,\dots, N_{last}\}$. $N_{last}$ is 7 unless otherwise specified.

Assignment constraints:
\begin{subequations}
\label{eq: app Scheduling assignment}
\begin{align}
     & \text{s.t.} 
     &&  \sum_{i \in I_R(r)} {{Ws}}_{i,n} \leq 1, \quad \sum_{i \in I_R(r)} {{Wf}}_{i,n} \leq 1, \quad r \in R, n \in N \\
     &&& \sum_{i \in I_R(r), n' \in N, n' \leq n} ({Ws}_{i,n'} - {Wf}_{i,n'}) \leq 1, \quad r \in R, n \in N  \\
     &&& \sum_{n \in N} ({Ws}_{i,n} - {Wf}_{i,n}) = 0, \quad i \in I  \\
     &&& \sum_{ i' \in I_R(r)} {Ws}_{i', n} \geq {Wf}_{i,n}, \quad r \in R, i \in I, n \in N, n < N_{last}, i \in I_R(r) \\
     &&& {Wf}_{i,0} = 0, {Ws}_{i, N_{last}}=0, \quad i \in I 
\end{align}
\end{subequations}
Variables $Ws_{i,n}$ and $Wf_{i,n}$ denote if the processing of a task i, belonging to the set of tasks $I_R(r)$ that can be scheduled on a certain machine $r \in R=\{M1, M2\}$, starts or ends at a scheduling event $n$. 

Duration, finish time, and time-matching constraints:
\begin{subequations}
\label{eq: app Scheduling duration}
\begin{align}
 & TCCH_{i,n} = \sum_{i' \in I} \kappa_{i,i'} Y_{i,i',n+1}, \quad i \in I, n \in N \setminus \{N_{last}\} \label{eq: app sch meh2}\\
     & TCCH_{i,N_{last}}=0, \quad T_{N_{last}}=MS  \label{eq: app sch meh1}\\
     & D_{i,n} = \alpha_i {Ws}_{i,n} + \beta_i {Bs}_{i,n}, \quad i \in I, n \in N  \\
     & Tf_{i,n} \geq Ts_{i,n} + D_{i,n} + TCCH_{i,n} - M_{big} (1-{Ws}_{i,n}), \quad i \in I, n \in N \\
     & Tf_{i,n} - Tf_{i,n-1} \leq M_{big} {Ws}_{i,n}, \quad i \in I, n \in N \setminus \{ 0 \} \\
     & Tf_{i,n} - Tf_{i,n-1} \leq D_{i,n} + TCCH_{i,n}, \quad i \in I, n \in N \setminus \{ 0 \} \\
     & Ts_{i,n} = T_{n}, \quad i \in I, n \in N \\
     & Tf_{i,n-1} \leq T_{n} + M_{big} (1-{Wf}_{i,n}), \quad i \in I, n \in N \setminus \{ 0 \} \\
     & Tf_{i,n-1} \geq T_{n} - M_{big} (1-{Wf}_{i,n}), \quad i \in I, n \in N \setminus \{ 0 \} \\
     & T_{n+1} \geq T_{n}, \quad n \in N \setminus \{N_{last}\}\\
     &  T_{0} = 0 
\end{align}
\end{subequations}
Variable $TCCH_{i,n}$ denotes the changeover duration for a task i at the end of the scheduling event n and depends on the sum of fixed changeover times $\kappa_{i,i'}$ from one task i to another $i'$. The changeover durations $TCCH_{i,n}$ and corresponding constraints were not initially accounted for in \citet{scheduling}. Variable $D_{i,n}$ is the duration of the processing of task i that starts at event n. $\alpha_i$ denotes the fixed set-up duration for a task i that starts at event n. $\beta_i$ denotes the variable duration per batch amount. There are three continuous variables that correspond to the batch size of task i at time point n: $Bs_{i,n}$ when it starts, $Bp_{i,n}$ when it is being processed, and $Bf_{i,n}$ when it finishes. $Ts_{i,n}$ and $Tf_{i,n}$ are the start and finish time of the processing of task i that starts at n. $T_n$ denote the time of scheduling event n. $M_{big}$ is a big-M constant. 

Batch size constraints:
\begin{subequations}
\label{eq: app Batch duration + MB}
\begin{align}
     & Bmin_{i} {Ws}_{i,n} \leq Bs_{i,n} \leq Bmax_{i} {Ws}_{i,n}, \quad i \in I, n \in N \\
     & Bmin_{i} {Wf}_{i,n} \leq Bf_{i,n} \leq Bmax_{i} {Wf}_{i,n}, \quad i \in I, n \in N \\
     &  Bmin_{i} (\sum_{n' \in N, n' < n } {Ws}_{i,n'} - \sum_{n' \in N, n' \leq n } {Wf}_{i,n'}) \leq  Bp_{i,n}\\
     & Bp_{i,n} \leq Bmax_{i} (\sum_{n' \in N, n' < n } {Ws}_{i,n'} - \sum_{n' \in N, n' \leq n } {Wf}_{i,n'}) \\
     & Bs_{i,n-1} + Bp_{i,n-1} = Bp_{i,n} + Bf_{i,n}, \quad i \in I, n \in N \setminus \{0\} \\
     & Bcons_{i,s,n} = \xi_{i,s} Bs_{i,n}, \quad s \in S, n \in N, i \in I_{s,in}(s) \\
     & Bprod_{i,s,n} = \xi_{i,s} Bf_{i,n}, \quad s \in S, n \in N, i \in I_{s,out}(s) \\
     & Bcons_{i,s,n} \leq Bmax_{i} \xi_{i,s} {Ws}_{i,n}, \quad s \in S, n \in N, i \in I_{s,in}(s) \\
     & Bprod_{i,s,n} \leq Bmax_{i} \xi_{i,s} {Wf}_{i,n}, \quad s \in S, n \in N, i \in I_{s,out}(s) \\
     & S_{s,n} = S_{s,n-1} + \sum_{i \in I_{s,out}(s)} Bprod_{i,s,n} - \sum_{i \in I_{s,in}(s)} Bcons_{i,s,n}, \quad s \in S, n \in N \\
     & S_{s,N_{last}} \geq P_{s}, \quad s \in S \\
     &  S_{s, 0} = S0_{s} - \sum_{s, i \in I_{in}(s)} Bcons_{i,s,0}, \quad s \in S \\
     & Cost_{CCH} = \sum_{i \in I, i' \in I, n \in N} Y_{i,i',n} \kappa_{i,i'} \label{eq: app Batch cost}
\end{align}
\end{subequations}
$Bmin_{i}$ and $Bmax_{i,n}$ denote the minimum and maximum batch sizes for a given task i. $Bcons_{i,s,n}$ and $Bprod_{i,s,n}$ are the amount of consumption and production of material $s \in S=\{I23, P1,P2,P3,P4\}$ in task i at event n and depend on the mass balance coefficient $\xi_{i,s}$. $S_{s,n}$ and $S0_s$ are the available and initial inventory of material s. The planning targets $P_s$ for each material s are set by the planning problem (apart from I23, which is set to 0). $Cost_{CCH}$ denotes the optimal changeover cost used to correct the planning-level scheduling cost estimate.

Tightening constraints:
\begin{subequations}
\label{eq: app Scheduling tightening}
\begin{align}
     & \sum_{n \in N, i \in I_{R}(r)} (D_{i,n} + TCCH_{i,n}) \leq MS , \quad r \in R \\
     & \sum_{n' \in N, i \in I_{R}(r), n' \geq n} (D_{i,n'} + TCCH_{i,n}) \leq MS - T_{n}, \quad r \in R, n \in N  \\
    & \sum_{n' \in N, i \in I_{R}(r), n' < n} (\alpha_i {Wf}_{i,n'} + \beta_i Bf_{i,n'}+TCCH_{i,n'}) \leq T_{n}, \quad r \in R, n \in N  \\
    & Y_{i,i',n} \geq {Wf}_{i,n} + {Ws}_{i', n} - 1, \quad r \in R, n \in N, i \in I_R(r), i' \in I_R(r), \kappa_{i,i'} > 0, n < N_{last}  
\end{align}
\end{subequations}

\subsection{Control formulation} \label{sec: Appendix ctrl opt}
All productions $m \in \{I3,P1,P2,P3,P4\}$ are modelled according to the chemical reaction $A \rightarrow B \rightarrow C$. A separate optimal control formulation (\ref{eq: app Control}) is formulated for each batch to optimize a mix of optimal processing time $t_f$ and final energy consumption $E(1)$ (\ref{eq: app Control obj}). Constraints (\ref{eq: app Control ca}) to (\ref{eq: app Control cc}) give the differential equations that determine the concentration profiles $c_A(t), c_B(t), c_C(t)$ as a function of time. The rate of reaction of $A$ to $B$ is first-order in the state variable concentration of A as well as the cooling flowrate control variable $u(t)$ with production-specific reaction constant $k_{1,m}$. Similarly, the rate of reaction of $B$ to $C$ is first-order in the concentration of B but proportional to $u(t)^{0.8}$ instead with reaction constant $k_{2,m}$. (\ref{eq: app control state1}) determines minimal terminal concentrations of B and C that are proportional to the batch target $Bprod_m$ set by the optimal scheduling decisions and normalized by parameters $V_B$ and $V_C$. (\ref{eq: app control state2}) denote the three initial concentrations, and (\ref{eq: app control state3}) denote path constraints on the ranges of $c_A(t), c_B(t), c_C(t)$. (\ref{eq: app Control en}) determines the differential expression with respect to time for the energy consumption $E(t)$ to be proportional with respect to the cooling flowrate $u(t)$ and to $V^2$. Finally, (\ref{eq: app control path}) limits the range of the independent time variable $t$, the control $u(t)$, and the processing time $t_f$.

\begin{subequations}
\label{eq: app Control}
\begin{align}
     & \min. 
     && 2 t_f + 0.5  E(1) \label{eq: app Control obj} \\
     & \text{s.t.} 
     &&  \frac{\partial c_A(t)}{\partial t} = - t_f k_{1,m} u(t) c_A(t) \label{eq: app Control ca} \\
     &&& \frac{\partial c_B(t)}{\partial t} = t_f ( k_{1,m} u(t) c_A(t) - k_{2,m} u(t)^{0.8} c_B(t))\\
     &&& \frac{\partial c_C(t)}{\partial t} = k_{2,m} u(t)^{0.8} c_B(t)) \label{eq: app Control cc} \\
     &&& \frac{\partial E(t)}{\partial t} = t_f V^2 u(t)  \label{eq: app Control en}  \\
     &&& c_B(1) \geq \frac{Bprod_m}{V_B}, \quad c_C(1) \geq \frac{Bprod_m}{V_C}  \label{eq: app control state1} \\
     &&& c_A(0) = 17, \quad c_B(0)= 0, \quad c_C(0)=0   \label{eq: app control state2}  \\
     &&& c_A(t), c_B(t), c_C(t) \in [0, 100]  \label{eq: app control state3}\\
     &&& t \in [0, 1], \quad u(t) \in [1, 9], \quad t_f \in [0, 500] \label{eq: app control path}
\end{align}
\end{subequations}

\section{Planning profiles for best solutions} \label{sec: Solution profiles}

%

\begin{figure}[htp!] 
    \centering
    \includegraphics[width=0.9\textwidth]{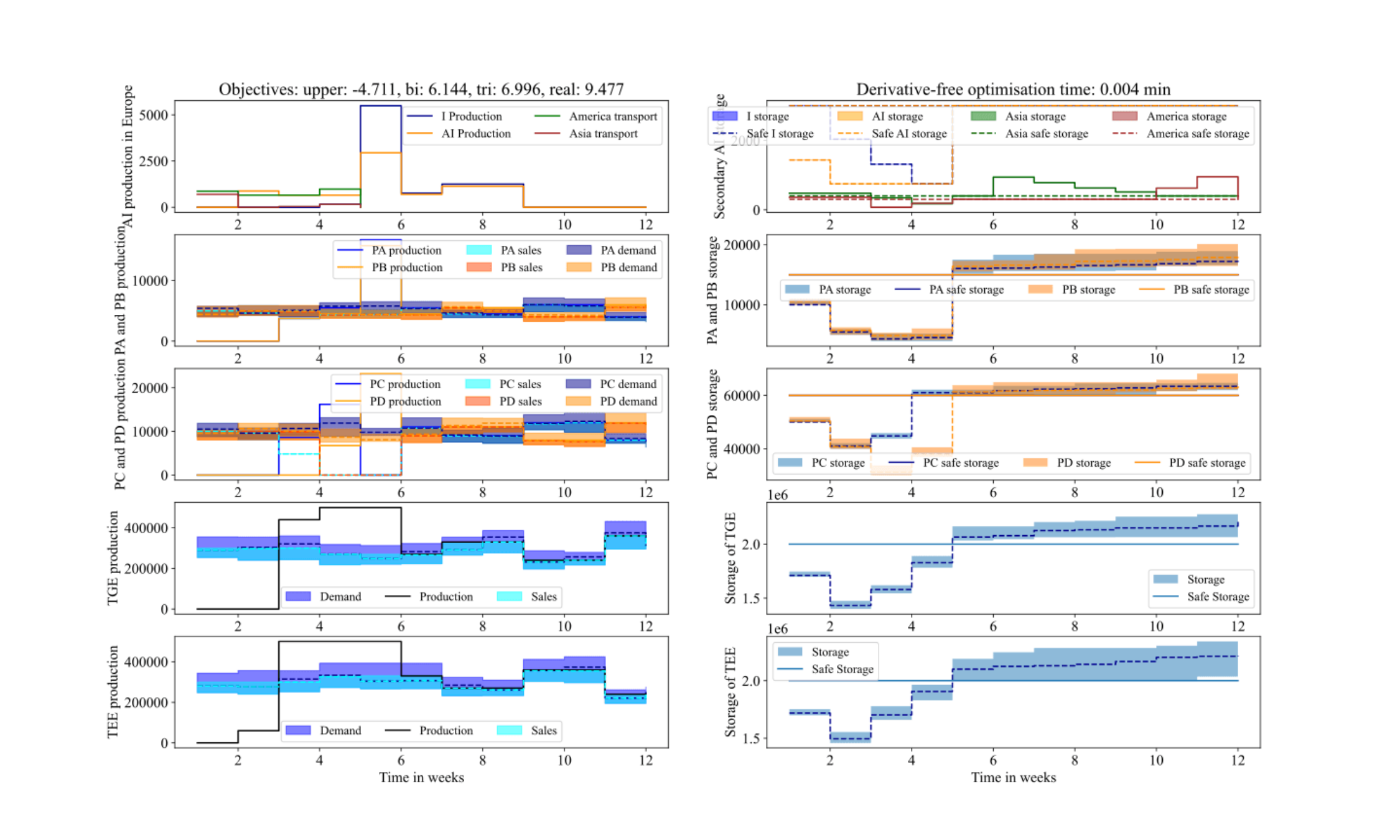}
    \caption{Solution dashboard and planning profile of the planning-only problem}
    \label{fig:profile_centralized}
\end{figure}

\begin{figure}[htp!] 
    \centering
    \begin{subfigure}[b]{0.9\textwidth}
        \centering
        \includegraphics[width=\textwidth]{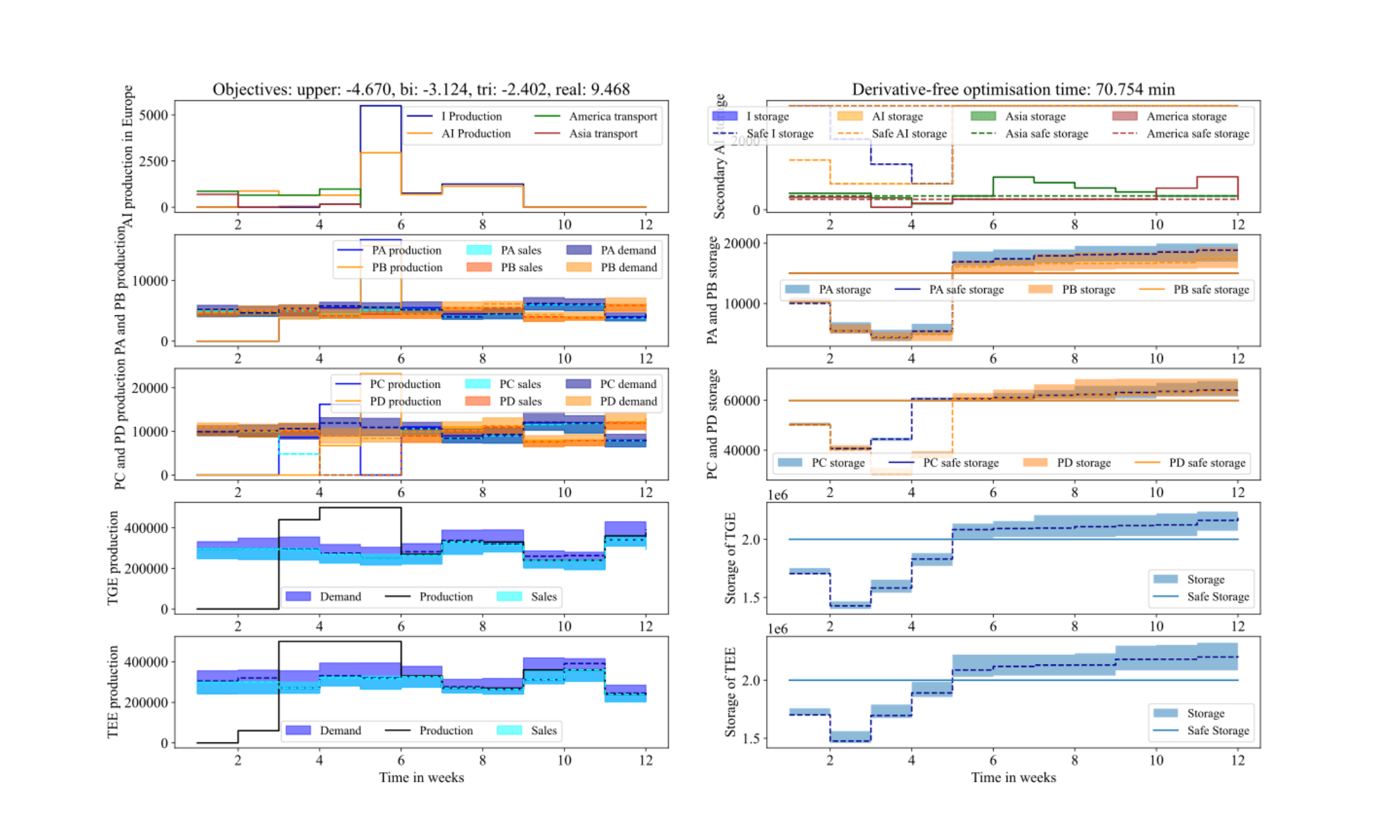}
        \caption{Solution dashboard and planning profile of the bi DFO approach on the      problem}
        \label{fig:profile_bi}
    \end{subfigure}
    \vskip \baselineskip
    \begin{subfigure}[b]{0.9\textwidth}
        \centering
        \includegraphics[width=\textwidth]{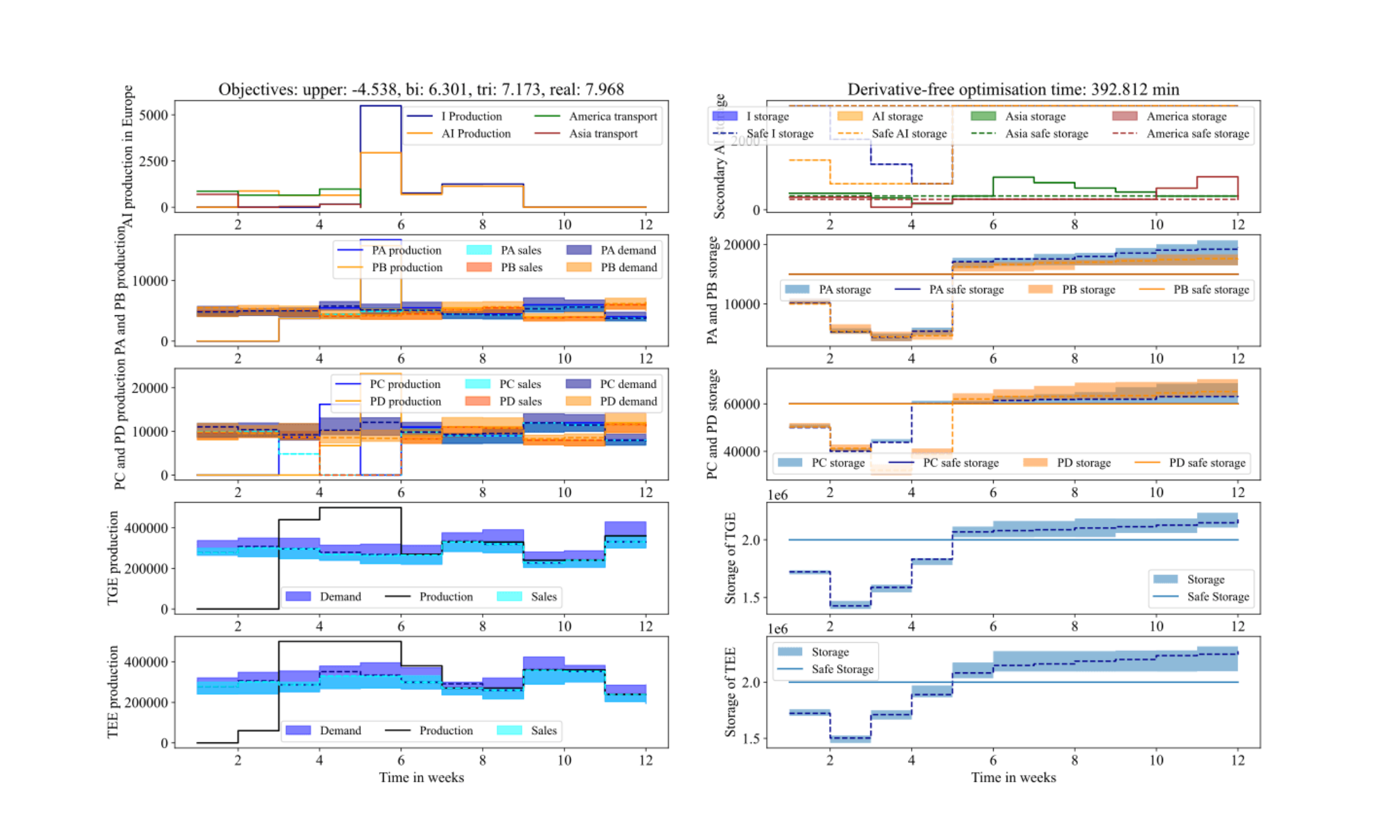}
        \caption{Solution dashboard and planning profile of the tri DFO approach starting   from the approx DFO solution}
        \label{fig:profile_bi_real}
    \end{subfigure}
    \label{fig:profiles1}
\end{figure}
\begin{figure}[htp!] 
    \begin{subfigure}[b]{0.9\textwidth}
        \centering
        \includegraphics[width=\textwidth]{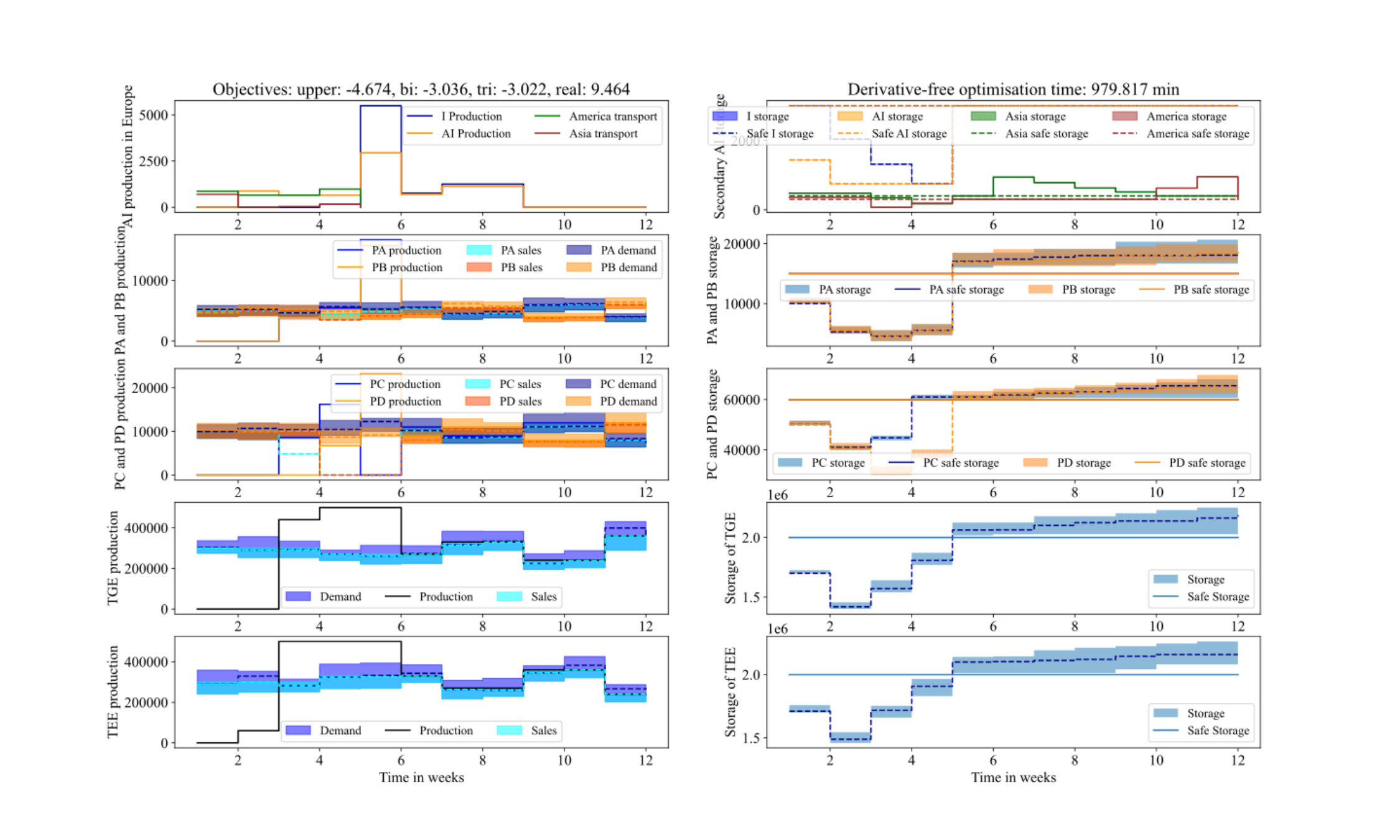}
        \caption{Solution dashboard and planning profile of the approx DFO approach on the      problem}
        \label{fig:profile_tri}
    \end{subfigure}
    \begin{subfigure}[b]{0.9\textwidth}
        \centering
        \includegraphics[width=\textwidth]{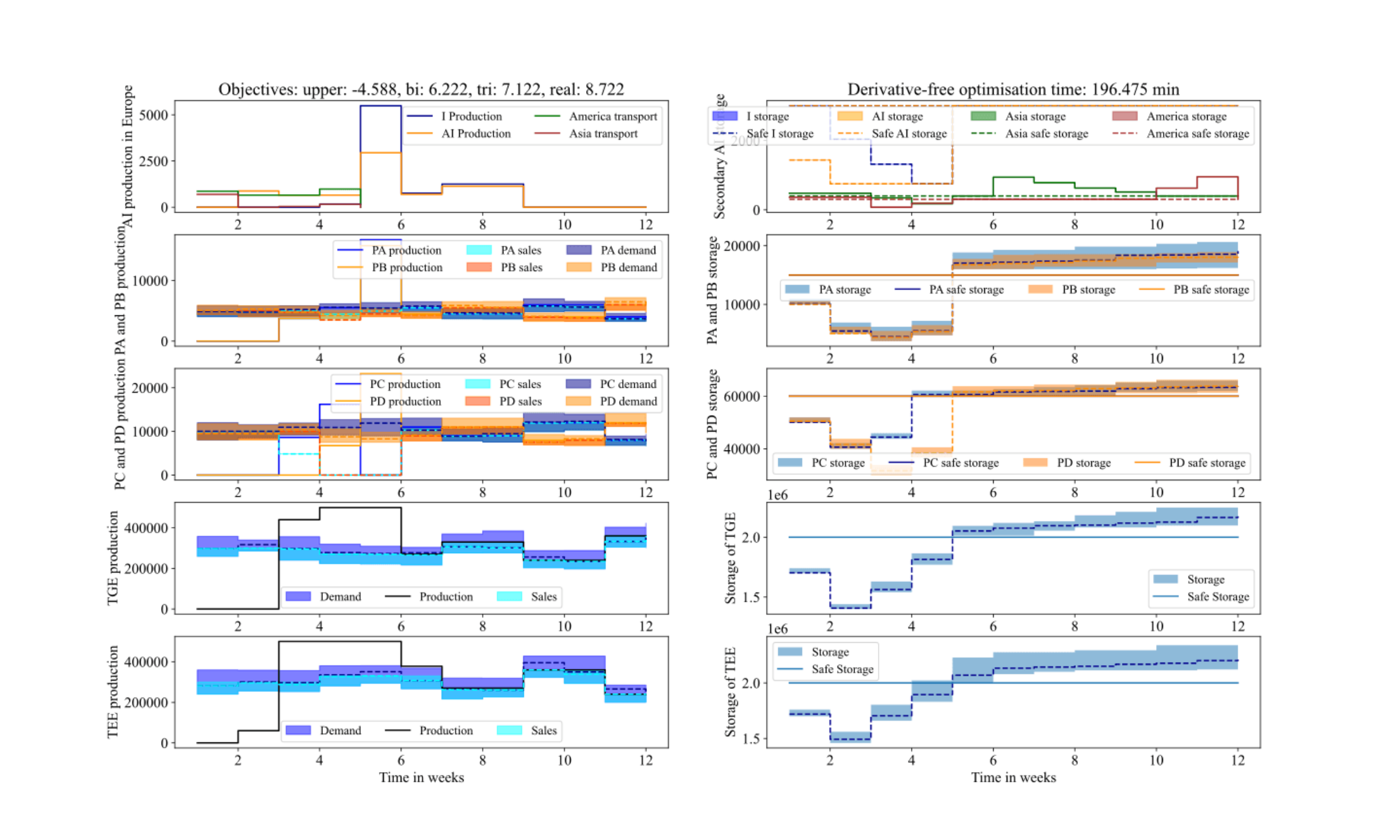}
        \caption{Solution dashboard and planning profile of the tri DFO approach starting   from the approx DFO solution}
        \label{fig:profile_real}
    \end{subfigure}
    \label{fig:profiles2}
\end{figure}

\begin{figure}[htp!] 
    \centering
    \begin{subfigure}[b]{0.9\textwidth}
        \centering
        \includegraphics[width=\textwidth]{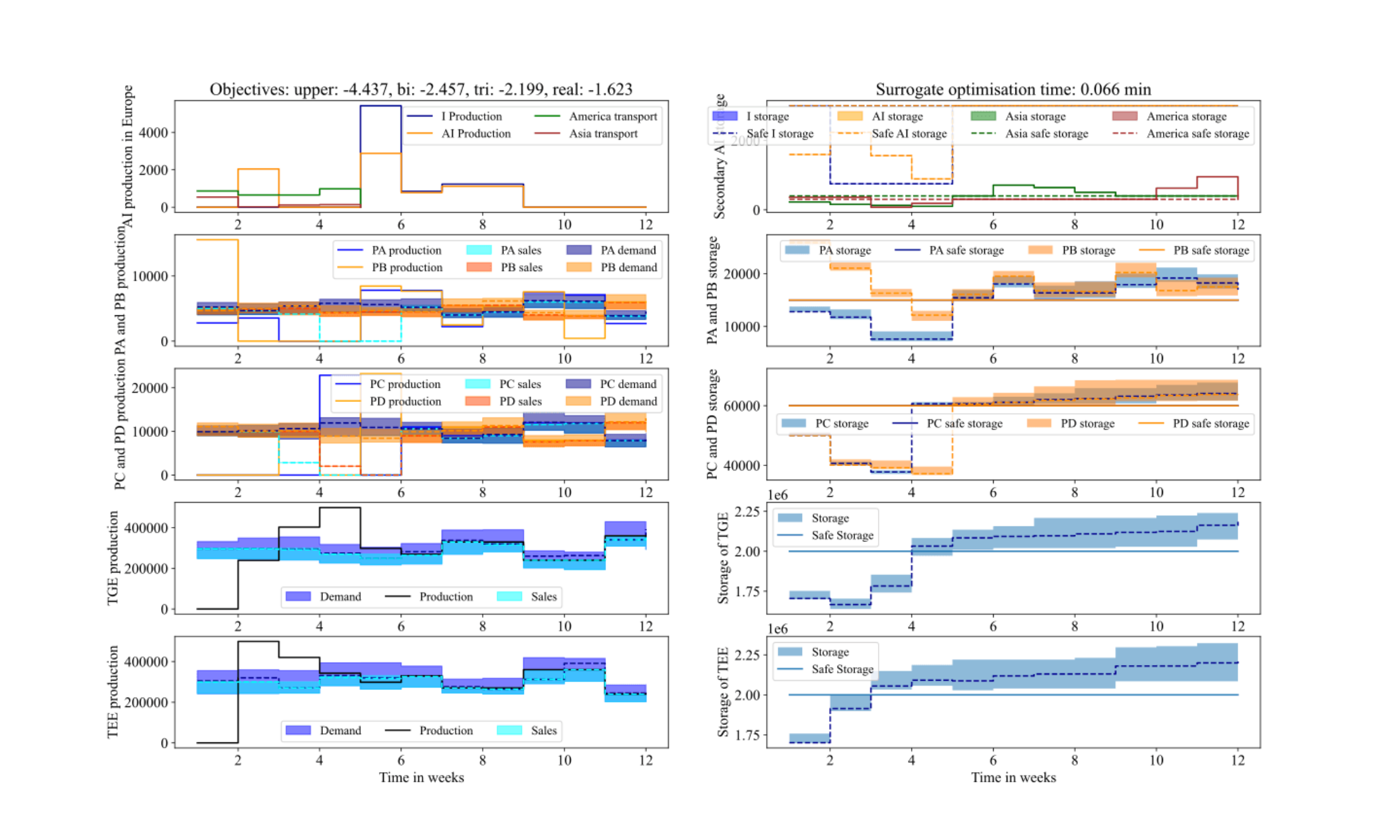}
        \caption{Solution dashboard and planning profile of the best approx surrogate   approach}
        \label{fig:profile_surr_integr}
    \end{subfigure}
    \vskip \baselineskip
    \begin{subfigure}[b]{0.9\textwidth}
        \centering
        \includegraphics[width=\textwidth]{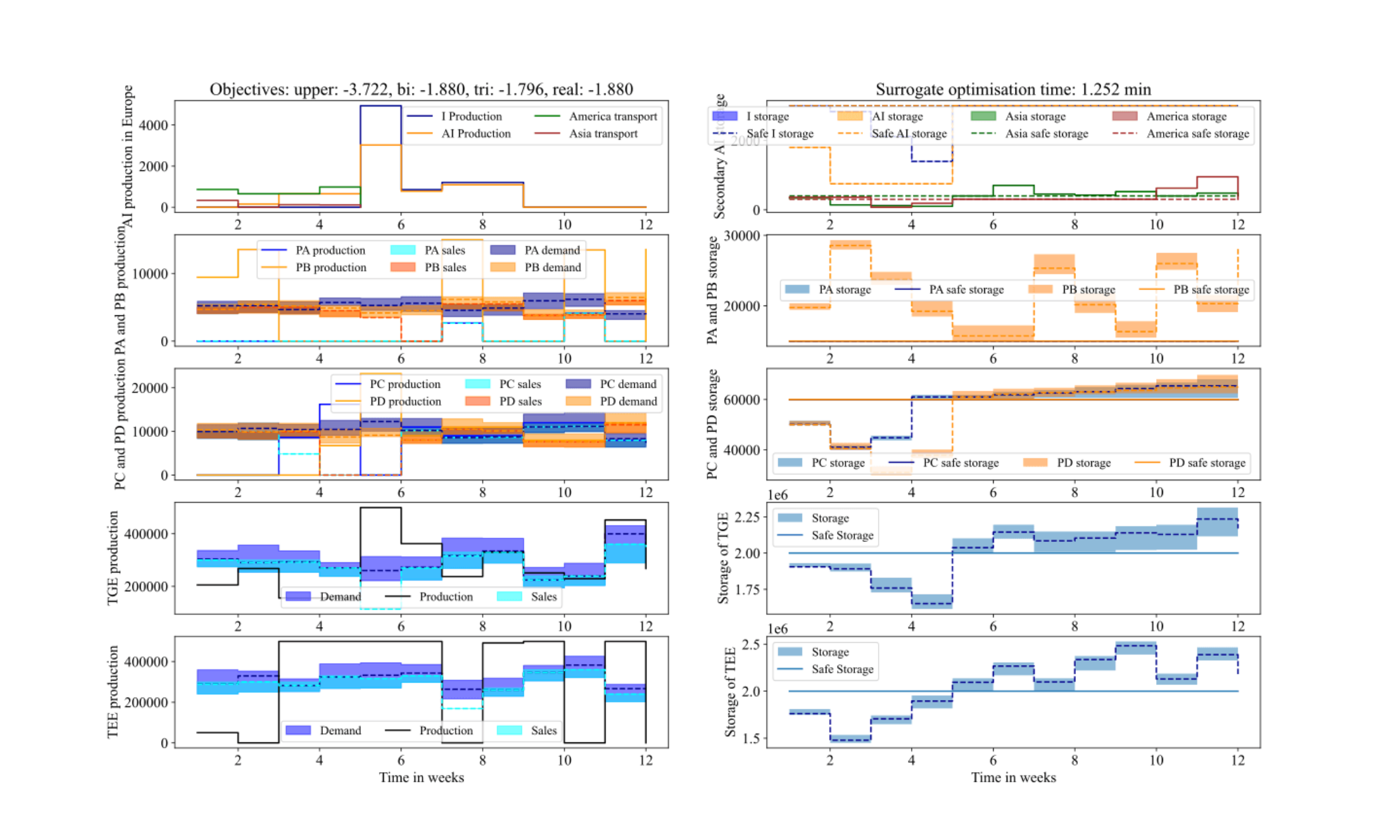}
        \caption{Solution dashboard and planning profile of the best tri surrogate approach}
        \label{fig:profile_surr_real}
    \end{subfigure}
    \label{fig:profiles3}
\end{figure}
\begin{figure}[htp!] 
    \begin{subfigure}[b]{0.9\textwidth}
        \centering
        \includegraphics[width=\textwidth]{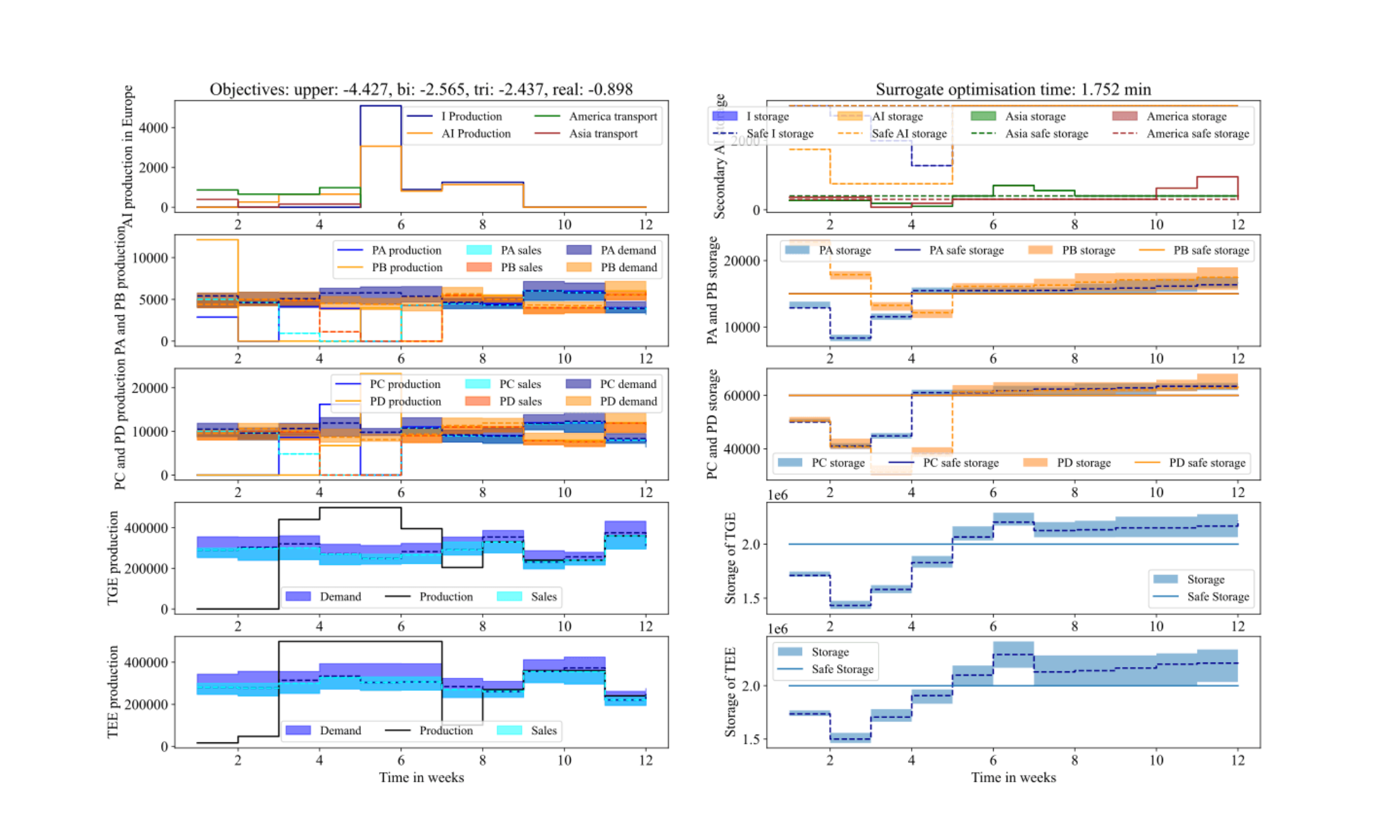}
        \caption{Solution dashboard and planning profile of the best bi surrogate approach}
        \label{fig:profile_surr_sch}
    \end{subfigure}
    \begin{subfigure}[b]{0.9\textwidth}
        \centering
        \includegraphics[width=\textwidth]{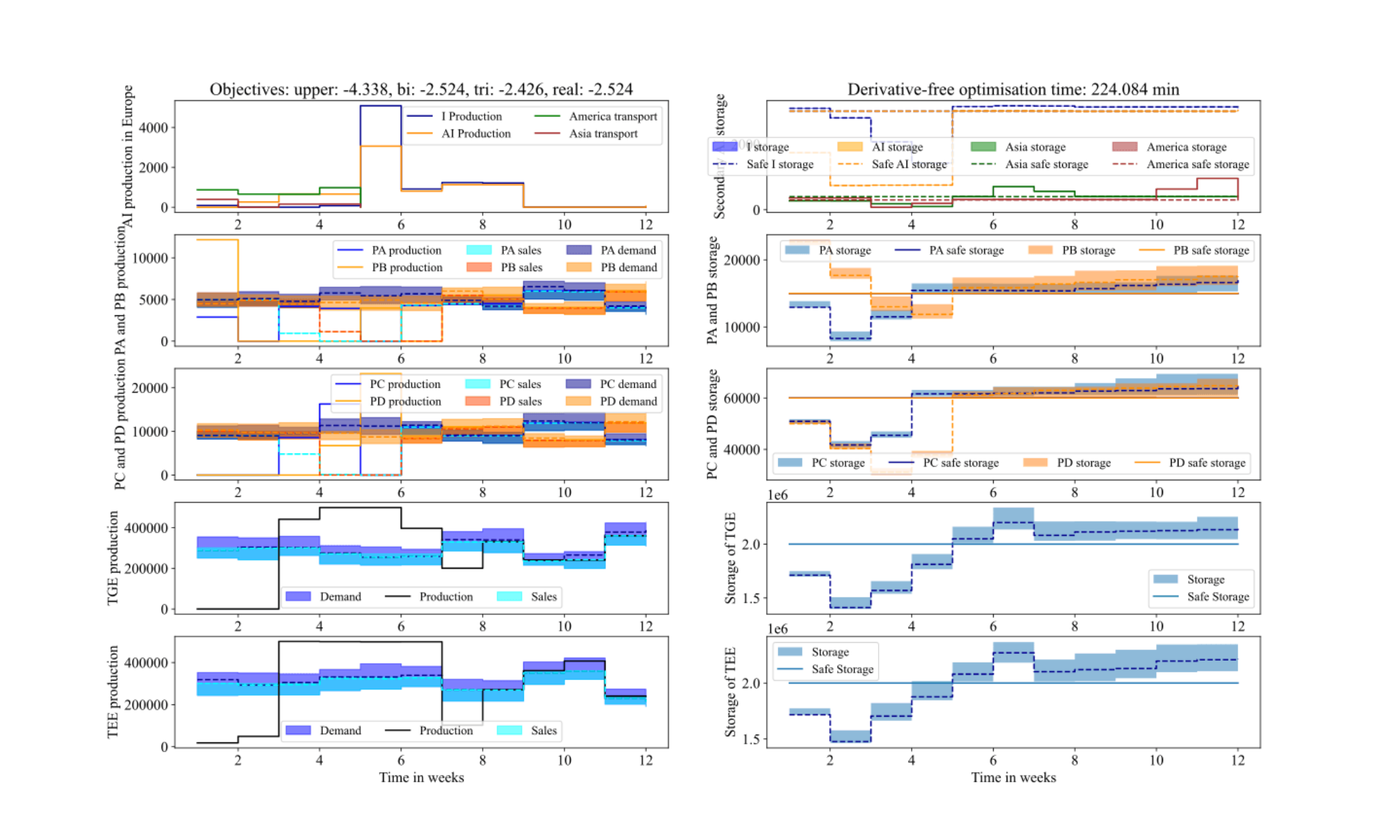}
        \caption{Solution dashboard and planning profile of the tri DFO approach starting   from the best bi surrogate approach}
        \label{fig:profile_hierarch_surr}
    \end{subfigure}
    \label{fig:profiles4}
\end{figure}

\end{document}